\documentclass{article}
\usepackage[utf8]{inputenc}

\usepackage{xr-hyper}
\usepackage{hyperref}
\usepackage{breakcites}
\usepackage{amsmath,amssymb,amsthm,pifont,bbm,csquotes}
\usepackage[dvipsnames]{xcolor}
\usepackage{booktabs}
\usepackage{adjustbox}
\usepackage{lscape,geometry,subcaption}
\usepackage{dsfont}
\usepackage{authblk}
\usepackage{enumerate}
\usepackage{etoolbox}
\usepackage{rotating}
\apptocmd{\thebibliography}{\setlength{\itemsep}{0pt}}{}{}

\usepackage[toc,page,titletoc]{appendix}
\usepackage[capitalise,nameinlink]{cleveref}

\hypersetup{
    colorlinks=true,
    linkcolor=blue,
    filecolor=magenta,
    urlcolor=cyan,
    pdftitle={Sharelatex Example},
    pdfpagemode=FullScreen,
}

\def \a    {\alpha}
\def \b    {\beta}
\def \d    {\delta}
\def \D    {\Delta}
\def \e    {\varepsilon}
\def \g    {\gamma}

\def \i    {\infty}
\def \l    {\lambda}

\def \oo {\"{o}}
\def \O    {\Omega}

\def \p    {\partial}
\def \s    {\sigma}

\def \grad {\nabla}

\newtheorem{thm}{Theorem}[section]
\newtheorem{lem}[thm]{Lemma}

\newtheorem{defn}[thm]{Definition}

\newtheorem{oppr}{Open Problem}

\theoremstyle{remark}

\newcommand{\eq}[1]{\begin{align*}#1\end{align*}}
\newcommand{\eql}[2]{\begin{align}\label{#2}#1\end{align}}
\newcommand{\as}[1]{\left\vert#1\right\vert}
\newcommand{\norm}[1]{\left\Vert#1\right\Vert}

\def \nn {\nonumber}

\usepackage[toc,page]{appendix}

\numberwithin{equation}{section}

\title{Open problems in PDE models for knowledge-based animal movement via nonlocal perception and cognitive mapping}

\author[1]{Hao Wang}
\author[1]{Yurij Salmaniw\footnote{The corresponding author. Email: salmaniw@ualberta.ca}}
\affil[1]{Department of Mathematical and Statistical Sciences, University of Alberta, Edmonton, AB T6G 2G1, Canada}
\date{}

\begin{document}
\maketitle
\vspace{-2ex}
\begin{abstract}
The inclusion of cognitive processes, such as perception, learning and memory, are inevitable in mechanistic animal movement modelling. Cognition is the unique feature that distinguishes animal movement from mere particle movement in chemistry or physics. Hence, it is essential to incorporate such knowledge-based processes into animal movement models. {Here}, we summarize popular deterministic mathematical models derived from first principles {that begin to incorporate such influences on movement behaviour mechanisms. Most generally, these models take the form of nonlocal reaction-diffusion-advection equations, where the nonlocality may appear in the spatial domain, the temporal domain, or both.} Mathematical rules of thumb {are} provided to judge the model rationality, {to aid in model development or interpretation, and to streamline an understanding of the range of difficulty in possible model conceptions}. {To emphasize the importance of biological conclusions drawn from these models, }we briefly present available mathematical techniques and introduce {some existing} ``measures of success" to compare and contrast the possible {predictions and} outcomes. {Throughout the review, we propose numerous open problems relevant to this relatively new area, ranging from precise technical mathematical challenges to broader conceptual challenges at the cross-section between mathematics and ecology. This review paper is expected to act as a synthesis of existing efforts while also pushing the boundaries of current modelling perspectives to better understand the influence of cognitive movement mechanisms on movement behaviours and space use outcomes.}
\end{abstract}
\tableofcontents

\section{Introduction}\label{sec:introduction}

The impact that the cognitive processes of organisms have on their movement is undeniable and ecologically important \cite{Fagan2013}. Cognitive processes, such as perception and memory, are unique features which distinguish the movement of animals from that of {non-living particles, such as atoms, molecules or projectiles}; similarly, the process of learning is a particularly unique feature which distinguishes merely directed movement, such as in anisotropic media \cite{Patlak1953} or aggregation of slime mould \cite{Keller1970}, from truly novel changes in behaviour and space use, a standard hallmark of learning \cite{ThrunPratt1998}. Without collecting and encoding information about the landscapes in which organisms live, many quintessential movement patterns, such as site fidelity and optimal foraging, would not manifest \cite{Fronhofer2013, van2009memory}. The apparent importance of cognition in animal movement processes warrants the development of mathematical models that incorporate these mechanisms \cite{smouse2010stochastic}. There already exist a number of reviews discussing the biological importance \cite{Fagan2013, Lewis2021}, the validity of such inclusions, as well as possible mechanisms behind the acquisition and use of available information \cite{Spiegel2016}. {Furthermore}, with the dramatic increase in types and quantities of animal movement data, {and} the significant decrease and increase in cost and computational power, respectively, various statistical methods and stochastic modelling efforts (such as \textit{individual based models} \cite{DeAngelis2005}) have grown in their use and application in this field \cite{Wilmers2015, Williams2020}. Such models will not be considered here, but not for their lack of importance or validity. Reviews of stochastic models with some focus on cognitive processes already exist \cite{smouse2010stochastic, Tang2010}, as well as at least one proposed ``standard protocol" to be used when considering such models in lieu of rigorous analytical techniques \cite{Grimm2006}. {Recent works such as \cite{Potts2022b} provide some interesting insights into the connections between individual-based models and partial differential equations, but we do not elaborate further in this work.} {To the best of our knowledge, this review is the first effort to {comprehensively} explore mathematical challenges {in the study of cognitive animal movement via partial differential equations (PDEs)}.}

{Here}, we aim to focus on deterministic models which incorporate some {forms} of cogniti{on}. We seek to address {both} \textit{how} to include certain processes from a mathematical standpoint, but also \textit{why} {these} formulations might correspond to their respective cognitive process. In doing so, we {hope to} encourage a balance between mathematical rigour and {maintenance} of realism. {Using the framework of PDEs, we can} incorporate and investigate the influence of explicit spatial structure on animal space use without appealing to a simpler ordinary differential equation structure, for example. On the other hand, a {PDE} setting is more {analytically} tractable than a stochastic or simulated setting, as {they offer} little means of analysis and therefore do not lend {themselves} to the discovery of ecological laws governing animal space use {outcomes}. This {potential} to perform analysis provides an additional layer of rigour through concrete and precise mathematical predictions, allowing one to answer some of the most important questions concerning animal movement and space use patterns. Indeed, explaining the spatial distribution of species using environmental factors has been named one of the top five ranked research fronts in ecology \cite{Renner2013}.

While a deterministic, continuous-in-time-and-space framework may be very difficult to validate and compare with {empirical} data, it is sometimes possible, see \cite{Moorcroft2,Ranc2020,NathanRanc2020, Ellison2020}. {Nevertheless}, even when these models cannot be fully integrated with empirical data, they still offer meaningful qualitative insights, as well as predictive power in the mechanisms {considered}. This is of particular importance in the area of {cognitive} movement ecology {as} we cannot (easily) directly observe what is happening in the brain of an organism. Even in cases where we \textit{can} observe certain brain functions or some other proxy \cite{Teitelbaum2016,Toledo2020,Fenton2020}, it is a more difficult challenge still to make the connection between these observations and the explicit behaviours of the organism. {At least one notable exception is the discovery of \text{place} and \textit{grid} cells \cite{Moser2015}, which are types of neurons that have been shown to be directly connected to external stimuli, such as landmarks or olfactory stimuli. {T}his provides a viable mechanism by which \textit{cognitive mapping} can occur (see Section \ref{sec:implicitmemory}).} {In most cases, though}, we can only make inferences on particular mechanisms of decision making based on the observed outcomes of the movements themselves. Hence, the validity of a proposed model or mechanism may at best correspond to its ability to accurately predict more general, qualitative trends in animal space use as observed in the available {empirical} data. While such comparisons may be lacking in precision, they can still provide meaningful insights and yield substantial motivation for future directions of research.

Consequently, our goals are the following. First, we will introduce some of the existing models within the framework of reaction-diffusion-advection equations. This will provide valuable context for less familiar readers, {and} motivation for more familiar researchers looking to extend these models in a meaningful way. {W}e hope to provide a reasonable amount of detail into the motivation behind the inclusion of certain modelling aspects, and how they connect with {natural phenomena in ecosystems} in an intuitive way. Second, we will discuss some of the predictions made and insights gained (in a biological sense) from each model {through mathematical analysis}. This closes the figurative loop through a connection between the mathematical constructions and the biological implications of each. {Throughout}, we will include generalizations and new formulations along the way. {We} compare and contrast {these new and} existing models so as to motivate and provide scaffolding to future researchers for further exploration of this exciting and still growing area of research, biologically and mathematically.

The remainder of this paper is organized as follows. In Section \ref{sec:PopModels} we introduce some of the existing prototypical movement models {featuring} forms of perception, memory and learning. {In Section \ref{sec:derivations}, we provide some useful mechanistic derivations for a general diffusion-advection equation along with the {corresponding} space-use coefficients.} In Section \ref{sec:biologicalinsights},  we discuss some important rules of thumb any modeller should consider, as well as the biological insights we can gain through study of the models introduced in Section \ref{sec:PopModels}. This includes a discussion of possible measures one may use to {explore parallels and discrepancies in predicted} space use outcomes. {We compliment the possible biological insights with more technical mathematical perspectives in Section \ref{sec:mathtech}.} To guide researchers moving forward, we provide many open {problems throughout the review}. We finish with {a broad summary of the key ideas presented and} some concluding remarks in Section \ref{sec:concremarks}.

\begin{figure}[ht]
\centering
\includegraphics[width=0.8\textwidth,height=10cm]{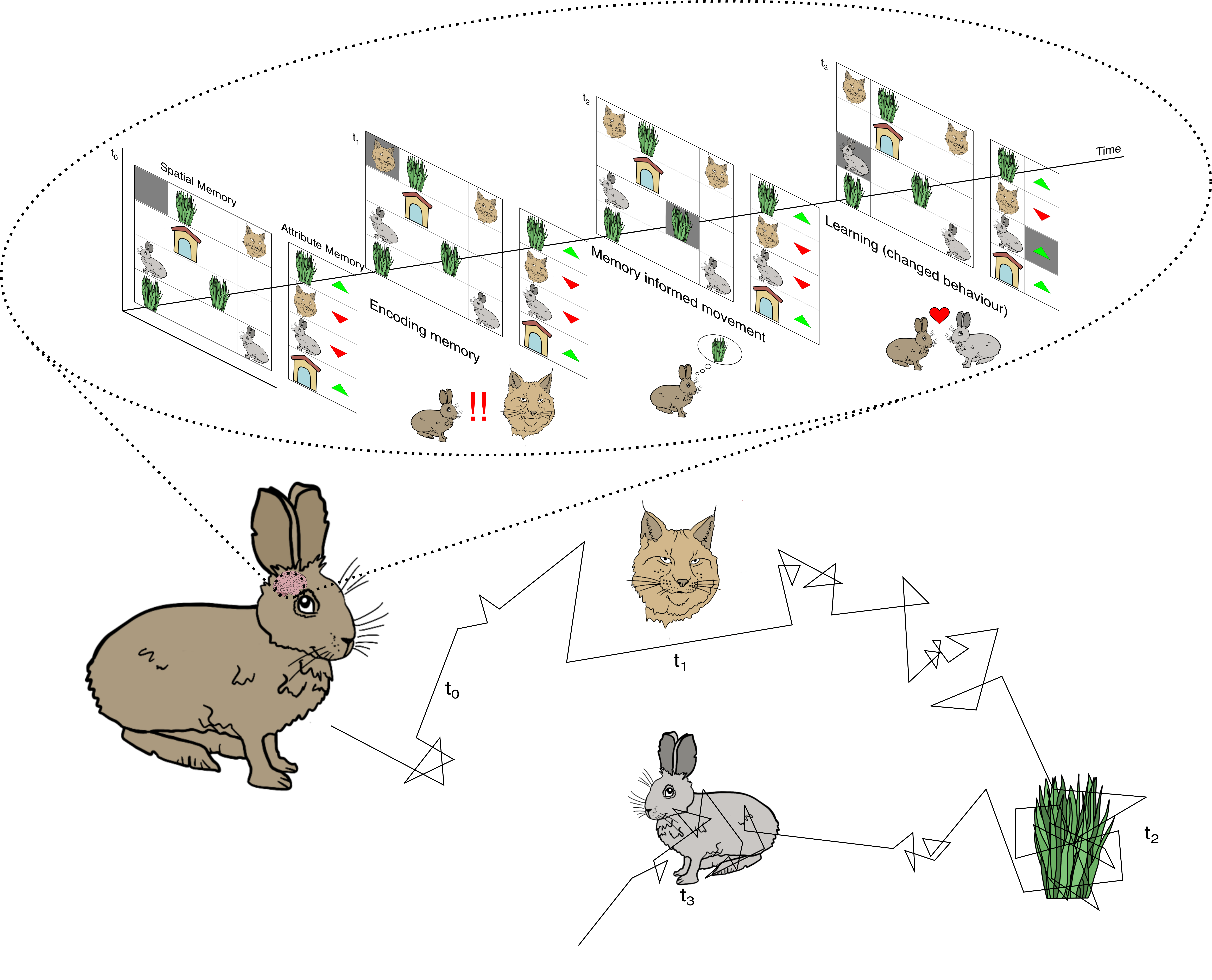}
\caption{Cognition in Animal Movement Conceptual Diagram}\label{ConDia}
\end{figure}

\section{Cognition in Animal Movement Models}\label{sec:PopModels}

{Before we investigate the precise mathematical formulation of cognitive animal movement models, we first conceptualize some of the key components and terminologies most useful in this context. There are three main categories of cognitive components commonly included in existing models: perception, memory, and learning. We appeal to the following definitions of each.}
{
\begin{defn}[Perception]
The ability to see, hear, or otherwise become aware of something through the senses.
\end{defn}
\begin{defn}[Memory]
The storage, retention and retrieval of information \cite{Fagan2013, Lewis2021}.
\end{defn}
\begin{defn}[Learning (psychology-based version)]\label{learning-psych}
The cause-effect process leading to information acquisition that occurs as a result of an individual's experience \cite{Lewis2021}.
\end{defn}
\begin{defn}[Learning (task-based version)]\label{learning-task}
The improved performance for a specific task as a result of prior experience \cite{Lewis2021}.
\end{defn}
} 

{
The first definition is certainly the least controversial; The note of caution, perhaps, is in the wide range of \textit{forms} of perception: visual and auditory cues are intuitively understandable for many, but numerous species have less familiar forms of perception, such as the ability to detect electromagnetic fields, chemical signals,  different wavelengths of light, or polarized light. Such considerations are important from the modelling perspective, which we explore further in Section \ref{sec:perception}.}
 
{
 The definition of memory is also widely agreed upon, though this may be a consequence of its broad scope; it does not make explicit the numerous subcategories of memory we know or suspect to exist. For this reason, we introduce two additional typologies, namely \textit{spatial} memory and \textit{attribute} memory.
 \begin{defn}[Spatial memory]
The memory for where objects/resources/places are in space \cite{Lewis2021}. Encodes spatial relationships or configurations \cite{Fagan2013, Lewis2021}.
 \end{defn}
  \begin{defn}[Attribute memory]
Encodes attributes of local features \cite{Fagan2013}.
 \end{defn}
}
{
Due to the complexity of memory and the various typologies that exist \cite{Schacter1992}, we focus our attention on these two forms of memory that are particularly pertinent to influencing animal movement. The distinction between these two forms of memory can be made clearer through a simple example: spatial memory encodes \textit{where} food is located, whereas attribute memory encodes the \textit{quality} and \textit{quantity} of the food. Of course, these two types often interact with each other, such as the storage of attribute information within a spatial structure. We expand on this idea in Section \ref{sec:implicitmemory}.}

{
The concept of learning is much less well defined, with definitions appearing to be discipline-dependent{; hence, w}e adhere to the most useful definitions in the context of informed animal movement From Definition \ref{learning-psych}, many of the models to be introduced here feature learning {implicitly}. From Definition \ref{learning-task}, we are more restricted in what may be considered learning{:} the learning process must be more explicit in some way, with measurable differences in outcomes a prerequisite. We expand more on these ideas in Section \ref{sec:learning}.}

{
Due to the importance of understanding the interplay {between} each of these components, we appeal to a conceptual diagram to {reinforce} these ideas. Figure \ref{ConDia} depicts {an} interplay between perception, spatial \& attribute memory, and learning in a home-range bound herbivore, and their influence on the animals' movement behaviours. The physical location of the animal in space at a given time step is denoted by the grey coloured box. Animals incorporate information about their environment by exploring and moving within their home range. These stimuli may come in the form of (a) food, (b) predators, (c) conspecifics, or even (d) the centre of its home range. The animal's attraction to (or repulsion from) these stimuli is dictated by its attribute memory, which assigns a quality to each landmark or stimulus the animal remembers. In this example the classification is binary, denoted by either a green arrow (attraction) or a red arrow (repulsion). All of this information is stored within the animals \textit{cognitive map} (see Section \ref{sec:implicitmemory}), which is stored within the animal's brain and {biases} future movement decisions. From the psychology-based definition of learning, the animal is going through a learning process as it updates its cognitive map (e.g., from $t_0$ to $t_1$, the animal has learned the location of a predator); from the task-based definition of learning, the animal is going through a learning process as it changes an attribute quality from repulsion to attraction (e.g., the change in attribute memory from $t_2$ to $t_3$). It is this dynamic interaction between animal movement and experience, environment, and its own cognitive abilities that we seek to model mechanistically.
}

We {now} introduce some of the popular modelling efforts which include {at least one } form of cognition. The order in which these results are presented is intended to, as best we are able, start with a more simplistic viewpoint before moving towards more complex formulations. {This} increase in complexity is primarily mathematical, {but as we will soon see,} the complexity of the cognitive function(s) included naturally escalates contemporaneously.

{T}o provide a strong foundation, we consider the following {prototypical} scalar advection-diffusion equation under a symmetric dispersal kernel
 \eql{
 \frac{\p u}{\p t} (x,t) &=  d \D u(x,t)  -   \grad \cdot \left(  u(x,t) \grad a(x,t) \right) 
 }{prototype}
 with a more general form  derived in {Section} \ref{sec:FPDer}. The diffusion rate $d$ {corresponds} to transition probabilities due to random movement, while the advective {potential} $a(x,t)$ {corresponds} to the bias in movement based on information at a spatial location $x$ at time $t$. In our {setting}, th{is} advective potential is the most important quantity we consider as it is where virtually all cognitive processes are currently incorporated. Heuristically, $a(x,t)$ can be thought of as the attractivity of a point at a given time, and so we may incorporate varying forms of cognition by adjusting this bias in movement through reasonable biological considerations. From a mechanistic point of view, this can be seen most clearly through the derivation of \eqref{prototype} where the quantity $a(x,t)$ is obtained through an exponential distribution of pertinent environmental covariates (see the derivation of space use coefficients found in {Section} \ref{sec:spacecoeffder} and \cite{Lele2006} for further discussion of resource selection functions). This should include perception, memory, learning, combinations of each, and their relation to other (external) environmental factors. From a mathematical point of view, a negative {sign} corresponds to attraction (moving \textit{up} the gradient of $a(x,t)$), while a positive sign corresponds to repulsion (moving \textit{down} the gradient of $a(x,t)$). Notice that this naturally allows {for both} attraction towards favourable regions {and} repulsion away from unfavourable regions. In what follows, we essentially derive the mechanisms by which these factors can be included through modification of the weighting function as derived in {Section} \ref{sec:spacecoeffder}.

\subsection{Perception}\label{sec:perception}

We start our exploration with a scalar equation model which includes an animal's ability to gather information about its landscape via \textit{nonlocal} perception, {motivating} subsequent models. In {many} scenarios, {an} ability to perceive is assumed to be based on purely local information{, which} may be an appropriate assumption for describing the movement of cells, for example. These local {advection} models{, such as the well-known Keller-Segel models \cite{Hillen2009,Painter2018}}, provide some initial motivation for how knowledge-based movement models {may} be constructed; for larger organisms, however, perception should not be {so} limited {as} it is well known that nonlocal cues, such as visual, auditory, olfactory or chemosensory cues, play a vital role in informing animal movement \cite{Fagan2017,Fagan2013,Painter2018}. Furthermore, the clarity by which animals can detect these cues may not be uniform across varying distances, let alone across species, within species, or even within individuals \cite{Zollner2000}. Consequently, there is substantial motivation to include nonlocal perceptual capabilities, and this should incorporate both \textit{distance} and \textit{quality of detection}. 

In {our} first example, the ``something" being detected is an assumed ``true" resource density function $m(x,t)$. It is assumed that the organism has a finite \textit{perceptual range} or \textit{detection scale} (the maximum distance at which landscape elements can be identified), as well as some description of how their perception {changes} with distance. Mathematically, this can be described by integration over space (a spatial convolution)
\eql{
h(x,t) := \int_\O m(y,t) g(x-y) dy,
}{resourcedetection}
which is sometimes referred to as a \textit{resource perception function} \cite{Fagan2017}. {W}e will refer to it more generally as a \textit{perception function}. {T}he kernel $g(x-y)$ describes the modifications in the forager's perception with distance, which we refer to as the \textit{perceptual kernel} or \textit{detection function}. {To clarify this distinction, we note that the perception function relates to \textit{what} the animal perceives, whereas the detection function relates to \textit{how} the animal perceives.} In \cite{Fagan2017}, the authors consider an unbounded, one dimensional spatial domain $\O = \mathbb{R}$ with the following possible perceptual kernels:
\eql{
g(x-y) &: =
\begin{cases} 
\frac{1}{2R} , \hspace{0.5cm}- R \leq x-y \leq R, \cr
0, \hspace{0.5cm} \text{otherwise},
\end{cases} \nonumber \\
g(x-y) &:= \frac{1}{\sqrt{2 \pi} R} e^{- (x-y)^2 / 2 R^2} , \\
g(x-y) &:= \frac{1}{2R} e^{- \as{x-y}/ R} \nonumber.
}{detectionkernels}
The quantity $R\geq0$ is the perceptual range, which is {proportional to} the standard deviation of the forager's detection function. These particular forms were chosen since the authors were interested in the transition between platykurtic (no tails) and leptokurtic (fat tailed) detection functions, each of which can be obtained from the exponential power distribution \cite{Smith1975}. In the first case, the so-called \textit{top-hat} detection function, {an} agent can perceive resources equally a fixed distance away from its current location and cannot {detect} beyond that fixed distance. The subsequent functions, the \textit{Gaussian} and \textit{exponential} detection functions, respectively, allow the agent to perceive nearby resources most clearly and decays monotonically as the distance from the observation location increases. In practice, a perceptual kernel could be any function satisfying the following {hypotheses}:
\begin{enumerate}[\hspace{0.1cm}i)]
\item $g ( x)$ is symmetric about the origin;
\item $\int_\O g (x) dx = 1$;
\item $\lim_{R \to 0^+} g (x) = \d (x)$;
\item $g(x)$ is non-increasing from the origin.
\end{enumerate}\label{detectionfunctions}
Condition i) assumes that the animal will perceive features equally across all directions. Condition ii) ensures that{,} given a perceptual range $R$, the mean perceptual range is the same for each detection function. Condition iii) assumes that as the perceptual range $R$ becomes arbitrarily small, the only information collected is purely local {(here, $\d(x)$ denotes Dirac's delta function{; for example, $h(x,t)$ converges to $m(x,t)$ as $R \to 0^+$})}. Finally, condition iv) assumes that an animals' perception does not improve as distance from the stimulus increases. This condition is mathematically convenient {while} also biologically reasonable, {but} it may be worth noting that some scenarios {do not satisfy} condition iv){, such as} \textit{hyperopia}, commonly referred to as farsightedness. All kernels introduced in \eqref{detectionkernels} satisfy these conditions. {Each kernel} can be generalized to any spatial dimension by turning $\as{x-y}$ into its corresponding norm in higher dimensional Euclidean space.

Starting from the prototypical model \eqref{prototype}, we may use the heuristic of $a(x,t)$ being the attractivity of a point so that $a(x,t) = \g h(x,t)$ itself becomes the attractive potential. {W}e reserve $\g \in \mathbb{R}$ to denote the strength of attraction to the potential $a(x,t)$. {In this case, we fix $\g$ positive since we are assuming that foragers will want to move to areas of higher resource density; in principle, $\g$ can certainly be negative, in which case foragers would be directed \textit{away} from high density areas.} The model with {nonlocal information gathering} and {exploratory movement} is then given by
\eql{
\frac{ \p u}{\p t} &= d \frac{ \p^2 u}{\p x^2} - \g \frac{ \p}{\p x} \left( u \frac{\p h}{\p x} (x,t)\right),
}{foragingsuccess}
 subject to the condition that $\int_\O u(x,t)dx = 1$. {T}his ensures the total population remains fixed, which is reasonable {as the model describes only movement}. In fact, this conservation is a consequence of the correct choice in boundary condition. A more general treatment of common boundary conditions and their implications are discussed in further detail in Section \ref{sec:considerations}. 
 
 In \cite{Johnston2019}, a similar {model} is obtained via a moment closure method to obtain a drift-anisotropic diffusion equation with focus on the 2-dimensional spatial case. This is done from the perspective of a velocity-jump random walk \cite{Othmer1988}, sometimes called a ``run-and-tumble" model, where an individual's movement is determined via a sequence of ``run" phases and ``turning" phases. The authors compare local and nonlocal gradient sampling with exclusive focus given to a uniform sampling over a circular region of radius $R$, which is exactly the top-hat detection function found in line \eqref{detectionkernels}. Their equation describing the motion of agents is identical to the general form found in {Section} \ref{sec:FPDer}. Similar to \cite{Fagan2017}, the authors then investigate the impact of local vs. nonlocal sampling under certain model types. 
 
 These works motivate a more general consideration of nonlocal detection and its influence on animal movement. Given a{n arbitrary} potential $a(x,t)$ {defined} in a domain $\O$, {the} perception function {is given} by
 \eql{
 \overline{a}_{g,R} (x,t) := \int_\O a(y,t) g(x-y) dy,
 }{generaldetection}
 where $g(\cdot)$ is any perceptual kernel satisfying the four {hypotheses} introduced {previously}.  We use the subscripts $g$ and $R$ to denote the dependence on the choice of perceptual kernel and perceptual radius, respectively.  {Notice that {to remain well defined} in a bounded domain, further modification may be required. We address such technicalities in more detail in Section \ref{sec:RoT}.}

 
 \subsection{Implicit Memory}\label{sec:implicitmemory}
 
 {W}e {now} explore how one may incorporate a rudimentary form of \textit{memory}. Memory plays a crucial role in the study of animal movement, yet remains a challenging problem for both biologists and mathematicians, as memory itself is a complicated process. Memories can be obtained via genetics (e.g., genetic triggers for migration, or inherited avoidance of predators \cite{Fagan2017}), or through direct experience. In this sense, memory is a higher level brain function than perception, as memory involves a secondary process of storing this {observed} information. We refer interested readers to \cite{Fagan2013} for an extensive review of the connections between memory and animal movement. Here, we seek to provide only key details {most} applicable to our modelling {efforts}. 
 
{As previously {noted}, we focus on spatial and attribute memory.} Of course, these two types often interact with each other, such as the storage of attribute information within a spatial structure. This process is sometimes referred to as \textit{cognitive mapping}. Originally, there was debate on whether or not cognitive maps exist; presently, the debate has shifted to what \textit{form} these maps actually take, (e.g., Euclidean vs topological, see \cite{Fagan2013}) and the references therein. Because we cannot (easily) directly observe these processes within the brain, models that include memory offer an {alternative} avenue to study these complicated agent-environment interactions. The challenge then becomes how to best model {a cognitive map}. {Most of the models to be introduced include spatial memory only; attribute memory is more difficult to include explicitly.} Indeed, it is more difficult to incorporate the \textit{quantity} of food in absolute terms since shifting a map by a large positive constant yields no difference in the model (i.e., the constant term vanishes since it appears underneath the gradient). {In this sense, attribute memory is included implicitly a priori since we make an initial assumption on whether a quantity is attractive or repulsive.} On the other hand, satisfaction measures discussed in Section \ref{sec:learning} may provide a useful avenue to study the effects of {explicit} attribute memory and {its} interplay with spatial memory. In the following sections, we discuss two differing perspectives that include a cognitive map, with {Section \ref{sec:staticmemory}} featuring static cognitive maps, and {Section \ref{sec:implicitdynamicmemory}} featuring dynamic cognitive maps. 

\subsubsection{Static Memory}\label{sec:staticmemory}
 
 The most obvious way to include memory is through a simple change of perspective in model \eqref{prototype}: define the quantity $a(x,t)$ to be the cognitive map of the animal. This could define desirable regions as well as regions to avoid, e.g., good resource locations or regions known to have predators. In this sense, model \eqref{foragingsuccess} could directly be interpreted as a \textit{memory} model if the function $m(x,t)$ is assumed to be the quantity being recalled. However, this approach may be viewed as naive as it requires the modeller to assume what the cognitive map actually looks like. In some cases, this may be more or less reasonable. For example, in order to study home range behaviour one may simply take $a(x,t) = a(x)$ to be the distance from a known home range site, such as a den, i.e., $a(x) = \g \Vert x - x_0 \Vert$, where $x_0$ is the fixed den site location, $\g > 0$ is the strength of attraction and $\Vert \cdot \Vert$ is the Euclidean norm. In this way, $\frac{\p a}{\p x}$ becomes a unit vector pointing towards the den site. This is precisely the formulation proposed in \cite{Moorcroft1}, which features an alternative derivation from a run-and-tumble perspective and an assumed Von-Mises distribution of turning directions; see also \cite{Johnston2019}. With a constant rate of diffusion, the model takes the form
 \eql{
  \frac{\p u}{\p t} &=  d \D u - \g   \grad \cdot \left( \frac{x - x_0}{\norm{x - x_0}} u \right).
 }{scalar-eqn-densite}
In this form the memory mechanism being included is {rather} rudimentary as it does not consider other factors that influence movement. However, since our general derivation in {Section} \ref{sec:spacecoeffder} includes the possibility of incorporating {several covariates}, a central den site may be considered {one} of many. 
 
 {A static cognitive map need not remain fixed in time; instead, it is static in the sense that it does not change based on the movement of the animal.} Given a resource density $m(x,t)$, it may be reasonable to assume that agents have some knowledge of the landscape relative to some other measure. Two such examples are the following: agents that are aware of the \textit{average} resource density, {and} agents that are aware of the \textit{per capita} resource density{, thereby assuming} that agents have knowledge beyond resource density {alone}. The first scenario has $a(x,t) = m(x,t)/\overline{m}$ {with movement} modelled by
 \eql{
 \frac{\p u}{\p t} = d \D u - \g \grad \cdot \left( u \grad \left( \int_\O \frac{m(y,t)}{\overline{m}} g(x-y) dy \right) \right) {= d \D u - \frac{\g}{\overline{m}} \grad \cdot \left( u \grad \overline{m}_{g,R} \right)},
 }{scalar-eqn-averagedensity}
 where $\overline{m} (t) = \frac{1}{\as{\O}} \int_\O m(y,t) dy$ is the average resource density at time $t$. The second scenario has $a(x,t) = m(x,t)/u(x,t)$ and can be modelled by
 \eql{
 \frac{\p u}{\p t} = d \D u - \g \grad \cdot \left( u \grad \left( \int_\O \frac{m(y,t)}{u(y,t)} g(x-y) dy \right) \right).
 }{scalar-eqn-percapitadensity}
 In \eqref{scalar-eqn-averagedensity}, we have a nonlocal equation that remains linear; \eqref{scalar-eqn-percapitadensity} is more complicated {as} it is nonlocal and nonlinear.
 
\begin{figure}[ht]
\centering
\includegraphics[width=\textwidth]{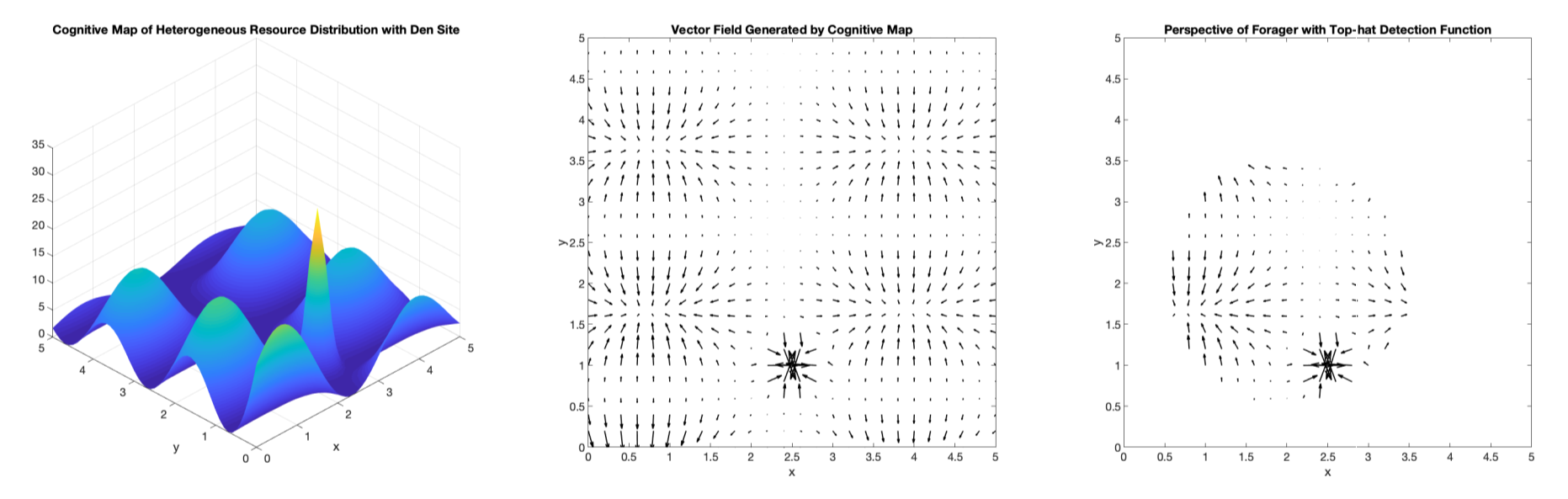}
\caption{A Sample Cognitive Map with Associated Vector Field and Perception}\label{StaticCogMap}
\medskip
\small
\end{figure}

Figure \ref{StaticCogMap} depicts a static, continuous-in-space cognitive map. The left panel is a sample static cognitive map with smaller peaks and troughs for high and low resources patches, in addition to a single tall peak denoting a den site. The middle panel provides the direction of movement based on the cognitive map, coming from the vector field generated by the advective potential. The right panel features the cognitive map with perception from the perspective of the forager at $(x,y) = (2,2)$ with {a top-hat detection function and} perceptual radius $R = 1.5$.
 
{A static cognitive map alone may not of significant interest; instead,} it {is more} interesting to consider the combination of a static cognitive map {with} a dynamic cognitive map. We consider other interacting cognitive maps through short and long term memory in Section \ref{sec:shortlongmemory}.

 \subsubsection{Dynamic Memory}\label{sec:implicitdynamicmemory}
 
 {In contrast} to a static cognitive map, one may consider a map to be a dynamically changing quantity, continuously updating as {an} agent moves throughout its environment. This offers more realism than the static cognitive map, as it is understood that memories are continuously formed and reformed as time passes. On the other hand, a dynamic cognitive map increases the mathematical complexity significantly as the description of movement for a single population may require {a second equation}.\\
 
  \noindent\textbf{Existing Models without Population Dynamics}
 
 The first form of a dynamic cognitive map we introduce is similar to the well-known group of \textit{Keller-Segel} models, which describe cell aggregation in response to chemical deposits left behind by the cells. {Rather than} following chemical deposits, it is {instead} assumed that the animals follow or avoid areas of high population density, assumed to be part of their cognitive mapping process. In general, we also include perception so that the {potential} $a(x,t)$ includes perception {of the population density}. The {cognitive map} is then $a(x,t) = \g \overline{u}_{g,R}$ for $\g \in \mathbb{R}$. The equation describing motion {of a single population} becomes
 \eql{
 \frac{\p u}{\p t} = d \D u - \g \grad \cdot \left( u \grad \overline{u}_{g,R} \right),
 }{scalar-eqn-1-1}
{where} $\g>0$ ($<0$) corresponds to attraction towards high-density {(low-density)} areas. Depending on the context, either could be valid: the first case may correspond to phenomena such as group defence strategies \cite{Shi_JDDE}, while the second may correspond to avoiding high-density areas where resources are expected to be less abundant.
 
 More interestingly, perhaps, are scenarios which include interactions between multiple populations. {To this end}, one may generalize model \eqref{scalar-eqn-1-1} to include $n$ interacting populations $u^i$, with the perception function for each group based on the varying population densities of all other groups:
 \eql{
 \frac{\p u^i}{\p t} = d_i \D u^i - \grad \cdot \left( u^i \grad \left(\sum_{j=1}^n \g_{ij} \overline{u}_{g,R}^j \right) \right).
 }{system-eqn-1}
 When the detection function $g(\cdot)$ is chosen to be the top-hat function, \eqref{system-eqn-1} is precisely the form proposed in scenario 1 of \cite{Potts2019} in a bounded, one-dimensional spatial domain. Similar to the single species model, the sign of $\g_{ij}$ determines whether species $u^i$ is attracted to or repelled from high population densities of population $u^j$. 
 
 Underlying such models is an implicit assumption that each population $u^i$ shares the same information, and so there must be some biological mechanism which allows agents to share this information between themselves. {Therefore}, this may be most applicable to very small organisms so that the density function $u^i$ is an appropriate description of how many individuals are found at a certain location in space. Recently, Potts {and Schl\"{a}gel} \cite{Potts2020} discussed this formulation in relation to step-selection analysis, as well as models' applicability and the potential for pattern formation to emerge. {We may now begin to formulate more precisely some open problems concerning model \eqref{system-eqn-1}.}
 
\begin{oppr} 
In what context do solutions exist solving the time-dependent problem \eqref{system-eqn-1} subject to periodic boundary conditions in a bounded domain? In \cite{Giunta2021}, a partial answer is established when the detection function $g(\cdot)$ is twice continuously differentiable (e.g., the Gaussian detection function), but few results exist for the top-hat detection function. The recent work \cite{Jungel2022} considers the existence of weak solutions to a nonlocal system of the form \eqref{system-eqn-1} and includes the case of the top-hat detection under the restriction that the nonlocal kernel components are in \textit{detailed balance} (see \cite[(H3)]{Jungel2022}. In the present context such a condition poses significant restriction, as this demands self-interaction (i.e., $\g_{ii} \neq 0)$) for every population, reducing potential for biological application. 
\end{oppr}

\begin{oppr}
It is clear that spatially constant steady state solutions exist solving problem \eqref{system-eqn-1}. In what context do spatially non-constant steady state solutions exist solving problem \eqref{system-eqn-1} when subject to zero-flux, homogeneous Neumann, or periodic boundary conditions?
\end{oppr}

\begin{oppr}
What are the qualitative properties of these non-constant steady states, whenever they exist? The recent work \cite{Giunta2022} provides some interesting insights and techniques in this direction.
\end{oppr}

\begin{oppr}
Do these solutions (either time-dependent or steady state) remain well-defined in the limit as $R \to 0^+$? In other words, can we meaningfully connect the nonlocal problem with its corresponding local problem in this limit? {\cite[Theorem 5]{Jungel2022} gives reason for optimism in the time-dependent case; on the other hand, the stability analysis of \cite{Potts2019} leaves less hope for the steady state problem due to the ill-posedness of the problem in this limit.}
\end{oppr}

\begin{oppr}
In what sense do solutions exist solving the time-dependent problem \eqref{system-eqn-1} subject to different boundary conditions? The challenge here is in appropriately defining the detection function near the boundary of the domain. See Section \ref{sec:considerations} for a more detailed discussion of boundary conditions.
\end{oppr}
 
 In models \eqref{scalar-eqn-1-1}-\eqref{system-eqn-1}, we have a dynamic cognitive map with perception without an additional equation. In this sense, these aggregation/segregation models are self-contained. In the following, the cognitive map is {now} described explicitly by an additional equation. This may be more appropriate to describe the cognitive map of larger organisms that have less dense populations, for example. The first example, {explored} in \cite{Lewis1993,Potts2019}, describes the movement of two or more populations in response to marks on the landscape (e.g., urine, faeces, footprints) left by the other population(s). To this end, denote by $p^i = p^i (x,t)$ the density of the presence of marks that are foreign to population $u^i$. It is assumed that marks grow linearly with respect to the presence of population $u^j$, and decay at a constant rate $\mu>0$. The evolution of marks foreign to population $u^i$ is given by
 \eq{
 \frac{\p p^i}{\p t} &= \sum_{j=1}^n \a_{ij} u^j - \mu p^i .
 }
Similar to previous examples, $\a_{ij} > 0$ indicates that population $u^i$ is attracted to population $j$, while $\a_{ij} < 0$ denotes repulsion ({note that \cite{Lewis1993,Potts2019} do not consider self-interaction, i.e., $\a_{ii} = 0$ is assumed}). Notice also that this form of ``memory" is somewhat different than our proposed definition, as their map is not stored within the brain, but {rather} within the environment itself. Despite this, it provides significant motivation for future models. 
 
 With this description of each $p^i$, we then {take} the attractive potential for population $u^i$ to be $a(x,t) = \g_i p^i(x,t)$ for some advection rate $\g_i \in \mathbb{R}$. {The} full model is given by
 \eql{
 \begin{cases}
  \frac{\p u^i}{\p t} = d_i \D u^i - \g_i \grad \cdot \left( u^i \grad  \overline{p}_{g,R}^i \right),  \cr
  \frac{\p p^i}{\p t} = \sum_{j=1}^n \a_{ij} u^j - \mu p^i  ,
  \end{cases}
 }{system-eqn-2}
 for $i = 1, \ldots, n$, where $n$ is the number of interacting populations. When the detection function $g(\cdot)$ is chosen to be the top-hat detection function, \eqref{system-eqn-2} is scenario 2 proposed in \cite{Potts2019}. Familiar readers may notice that if $n=1$, \eqref{system-eqn-2} is very similar to a Keller-Segel system in the limit as $R \to 0^+$. The key difference is the lack of diffusion appearing in the equation for $p^i$: diffusion, in this setting, has a regularizing effect, and so an absence of diffusion increases the difficulty of analysis (see Section \ref{sec:wellposedness} for further {discussion}).

 \begin{oppr}
In what sense do solutions exist solving the time-dependent problem \eqref{system-eqn-2} subject to periodic boundary conditions and a top-hat detection function in a bounded domain? Do solutions remain well defined in the limit as $R \to 0^+$?
 \end{oppr}

  \begin{oppr}
It is clear that spatially constant solutions exist as steady states of problem \eqref{system-eqn-2}. In what sense do spatially non-constant solutions exist solving the steady state problem associated to \eqref{system-eqn-2} subject to periodic boundary conditions and a top-hat detection function in a bounded domain? Do solutions remain well defined in the limit as $R \to 0^+$? 
 \end{oppr}

\begin{oppr}
What are the qualitative characteristics of the spatially non-constant steady state solutions obtained for problem \eqref{system-eqn-2}, whenever they exist? How do these qualitative properties compare to steady state solutions of model \eqref{system-eqn-1}? Of particular interest is the following general question: are the dynamics of model \eqref{system-eqn-2} included in the dynamics of model \eqref{system-eqn-1}? In other words, are new dynamics observed when the complexity is increased through explicit description of the cognitive map?
\end{oppr}

 \begin{oppr}
In what sense do solutions exist solving problem \eqref{system-eqn-2} subject to zero-flux, homogeneous Neumann, or homogeneous Dirichlet boundary conditions with a top-hat detection function in a bounded domain?
 \end{oppr}
 
 {Different from marks on a landscape,} in the following model the cognitive map {is recorded} within the foragers mind. To achieve this, the main idea is to track direct encounters between agents from different populations, the areas at which these occur referred to as \textit{conflict zones}. It is assumed that each population remembers an area where a conflict has occurred, and will be more likely to avoid this area in the future. Should they return to a location and experience no conflict, the cognitive map is updated accordingly. It is also assumed that memory decays at some rate proportional to the time since an event has occurred. This can be viewed as a combination of attribute memory and spatial memory, where the conflict is the attribute recorded at some spatial location where the conflict occurs.
 
 Denote by $k^i(x,t)$ the spatial memory of conflict zones for population $u^i(x,t)$. For simplicity, we first consider the case of two interacting populations. From our preliminary assumptions, $k^i$ should grow with respect to interactions between $u^1$ and $u^2$, while it should decay proportionally to $u^i$ and linearly some rate $\mu>0$. The equation describing the evolution of the spatial cognitive map then takes the form
 \eq{
 \frac{\p k^i}{\p t} &= \rho u^1 u^2 - ( \mu + \b u^i ) k^i .
 }
 The quantity $\rho \geq 0$ is the rate at which encounters occur; $\mu\geq 0$ is the rate at which memory decay with time; $\b\geq0$ is the rate at which the conflict zone decays due to agents revisiting the site and experiencing no conflict. 
 
 We now take a moment to note an important distinction between this model introduced above and the similar form introduced in \cite{Potts2016}. In the cases introduced in this review, we focus on describing the cognitive map as a \textit{magnitude}, describing important areas versus less important areas in relative terms. On the other hand, some works formulate the cognitive map as a \textit{probability density}, and so the specific form on the dynamic cognitive map is slightly different. In this case, the equation for $k^i$ above is instead derived to be
 \eql{
  \frac{\p k^i}{\p t} &= \rho u^1 u^2 (1-k^i) - ( \mu + \b u^i ) k^i .
 }{eqn1.123}
 In this way, outcomes are treated similar to a coin flip: a location is either part of a conflict zone, $(1-k^i)$, or it is not, $k^i$. While this is a subtle difference in interpretation, {it is not clear {whether}} the overall dynamics {should} appear roughly the same. Since the more popular method is to describe the cognitive map as a magnitude, we focus on these cases instead.
 
This can readily be generalized to $n$ interacting {populations} {as follows}:
 \eq{
 \frac{\p k^i}{\p t} &= u^i \sum_{j\neq i} \rho_{ij} u^j - ( \mu + \b u^i ) k^i ,
 }
 where $\rho_{ij}$ now denotes the rate at which encounters occur between populations $u^i$ and $u^j$. Then, for each $u^i$, the attractive potential is the cognitive map $k^i$ combined with perception, i.e., $a^i(x,t) = \g_i \overline{k}^i_{g,R} (x,t)$, where $\g_i \geq 0$ denotes the rate at which species $i$ moves away from all other populations (notice the sign change of the advective term below, recalling that positive advection denotes repulsion). The full model describing the evolution of $n$ interacting {populations} remembering conflict zones with perception is then given as
  \eql{
 \begin{cases}
  \frac{\p u^i}{\p t} = d_i \D u^i + \g_i \grad \cdot \left( u^i \grad  \overline{k}_{g,R}^i \right),  \cr
  \frac{\p k^i}{\p t} = u^i \sum_{j=1}^n \rho_{ij} u^j - (\mu + \b u^i) k^i  .
  \end{cases}
 }{system-eqn-3}
When the kernel $g(\cdot)$ is taken to be the top-hat detection function, this is scenario 3 proposed in \cite{Potts2019} in a bounded, one-dimensional spatial domain. 
 
\begin{oppr} 
In what sense do solutions exist solving the time-dependent problem \eqref{system-eqn-3} subject to periodic boundary conditions and top-hat function in a bounded domain? In what sense do spatially non-constant steady state solutions exist? What are the qualitative properties of these solutions? Do solutions remain well defined in the limit as $R \to 0^+$?
\end{oppr}

\begin{oppr}
In what sense do solutions exist solving the time-dependent problem \eqref{system-eqn-3} subject to zero-flux, homogeneous Neumann, or homogeneous Dirichlet boundary conditions with a top-hat detection in a bounded domain?
\end{oppr}

\begin{oppr} 
What are the key differences, if any, between solutions obtained for model \eqref{system-eqn-3} and solutions obtained for the same problem when the cognitive map is instead of the form described in \eqref{eqn1.123}?
\end{oppr}
 
 In the models presented so far, we have discussed the cognitive processes in movement population models that consider only movement (no birth or death of the population) between interacting populations. Sometimes, this is justifiable if one assumes that the movement process occur at a timescale that is much faster than that of a birth/death process. Readers should take caution, however, since such an assumption may invalidate the use of a \textit{quasi-steady state approximation}. We discuss this point in more detail in Section \ref{sec:RoT}. This observation motivates the following open problem.
 \begin{oppr}
 How do the dynamics of any model introduced in Section \ref{sec:implicitdynamicmemory} change when birth/death processes are also included? Do solutions exist solving the time-dependent or steady state problems? How do solution profiles change for spatially non-constant steady state solutions?
 \end{oppr}
 
\noindent\textbf{Existing Models with Population Dynamics} 
 
 Next, {we consider }a classical consumer-resource model with an additional term biasing the movement of the consumer. A slightly more general formulation is considered in \cite{Song2021Preprint} which is currently under review. Denote by $u(x,t)$, $v(x,t)$ the consumer and resource, respectively. The case {considered} is the most straightforward: {we assume} that the consumers have knowledge of where the resources are. We then take the perception function $a(x,t) = \g \overline{v}_{g,R} (x,t)$ for $\g > 0$. The consumer-resource model with knowledge and perception of resources is then described by
 \eql{
  \begin{cases}
 \frac{\p u}{\p t} = D_1 \D u - \g \grad \cdot \left( u \grad \overline{v}_{g,R} \right)  + \frac{c \b u v}{\a + v} - d u, \\
\frac{\p v}{\p t} = D_2 \D v + r v (1 - v/K ) - \frac{\b u v}{\a + v}.
 \end{cases}
 }{system-eqn-RC1-case1}
 This perspective is comparable to model \eqref{system-eqn-1}: instead of knowledge of the (current) density of other populations, the consumers  have knowledge of the current resource density. Here, $r>0$ denotes the maximum reproduction rate for the resource, while $K >0$ is the carrying capacity for the resource. The consumer is assumed to grow according to a Holling type II functional response with the growth rate $\b>0$ and the half-saturation constant $\a>0$, and decay linearly at rate $d>0$. The quantity $c > 0$ is a conversion efficiency by the consumer from the resource. Notice that if one takes $\g=0$, the system is reduced to a classical consumer-resource model. {On the other hand, in the limit as $R\to0^+$ the model reduces to a standard predator-prey model with prey-taxis, see \cite{Wang2021} and the references therein.}

 \begin{oppr}
In what sense do solutions exist solving the time-dependent problem \eqref{system-eqn-RC1-case1} subject to periodic boundary conditions in a bounded domain? Some insights can be found in \cite{Song2021Preprint}, however the authors considered a no-flux boundary condition in a bounded, one-dimensional spatial domain for mathematical analysis, and considered periodic boundary conditions in simulations.
 \end{oppr}
 
 \noindent\textbf{New Models and Extensions}
 
 In what follows, we discuss some new models that have not yet been investigated in the presented format. Some may require further development, {but} we still include them to motivate future research with the advancements currently being made. The first, moderately simple generalization is applied to model \eqref{system-eqn-3}. Some authors argue that a memory process should carry a similar derivation to the movement process itself \cite{Gourley2002}, which in this case suggests that the cognitive map should also feature some rate of diffusion. This leads to the idea of memory \textit{smearing}, which {introduces} random errors in memory recall. In this way, the diffusion smears the memory component $k^i$ so that {memory} is roughly accurate, but not remembered exactly{, and the imprecision increases over time until reinforced further}. This simple modification results in the new model
 \eql{
 \begin{cases}
   \frac{\p u^i}{\p t} = d_i \D u^i - \g_i \grad \cdot \left( u^i \grad  \overline{k}_{g,R}^i \right),  \cr
  \frac{\p k^i}{\p t} = \e_i \D k^i + u^i \sum_{j=1}^n \rho_{ij} u^j - (\mu + \b u^i) k^i ,
  \end{cases}
 }{system-eqn-3-mod1}
 where $0 < \e_i \ll d_i$. The {parameter} $\e_i$ is meant to include this mechanism of memory smearing, or how memories may be altered with respect to distance \cite{bracis2015memory}. 
 
\begin{oppr} 
In what sense do solutions exist solving the time-dependent problem \eqref{system-eqn-3-mod1} subject to periodic boundary conditions and top-hat function in a bounded domain? Note that since the cognitive map equation features a diffusive term, solutions are expected to be more regular than cases without diffusion.
\end{oppr}

\begin{oppr}
In what sense do solutions of the problem \eqref{system-eqn-3-mod1} converge to solutions of problem \eqref{system-eqn-3} in the limit as $\e_i \to 0^+$?
\end{oppr}
 
In model \eqref{system-eqn-3}, it is assumed that conflicts occur at exactly one point. However, the inclusion of perception may also be relevant in this term, since agents may be able to experience a conflict at a distance. This could result in a conflict (through a ``stand-off") that should also be remembered. Hence, the growth term for $k^1$ should feature the form $\rho \overline{u}_{g,R}^1 u^2$ in the case of two interacting species. For $n$ interacting species, this {map} takes the form
   \eql{
  \frac{\p k^i}{\p t} = \overline{u}^i_{g,R} \sum_{j=1}^n \rho_{ij} u^j - (\mu + \b u^i) k^i  .
 }{system-eqn-3-mod2}
 Naturally, these models assume the same perceptual kernel and perceptual radius for each population. This could be a reasonable assumption for studying more uniform populations, such as animals within the same species (e.g., wolf packs), but this may not be an accurate description if the interacting populations are significantly different.

\begin{oppr}
In what sense do solutions exist solving the time-dependent problem \eqref{system-eqn-3-mod2} subject to periodic boundary conditions and a top-hat detection function in a bounded domain?
\end{oppr}

\begin{oppr}
In what sense do spatially non-constant steady state solutions exist associated to problem \eqref{system-eqn-3-mod2} subject to periodic boundary conditions and a top-hat detection function in a bounded domain? What qualitative properties do these solutions hold?
\end{oppr}
 
\begin{oppr}  
How does nonlocal information gathering effect space use outcomes? More precisely, how do the qualitative features of spatially non-constant steady state solutions found for problem \eqref{system-eqn-3-mod2} differ from those found for problem \eqref{system-eqn-3}?
\end{oppr} 

\begin{oppr}  
In what sense do solutions exist solving the time-dependent problem \eqref{system-eqn-3-mod2} subject to zero-flux, homogeneous Neumann, or homogeneous Dirichlet boundary conditions with a top-hat detection function in a bounded domain?
\end{oppr} 

 Next, we discuss some extensions of the consumer-resource model prototype introduced in \eqref{system-eqn-RC1-case1}. In general, model \eqref{system-eqn-RC1-case1} can be written as
 \eql{
 \begin{cases}
 \frac{\p u}{\p t} = D_1 \D u - \g \grad \cdot \left( u \grad a(x,t) \right)  + \frac{c \b u v}{\a + v} - d u, \\
 \frac{\p v} {\p t} = D_2 \D v + r v (1 - v/K ) - \frac{\b u v}{\a + v},
 \end{cases}
 }{system-eqn-RC1}
 for some attractive potential $a(x,t)$. We can then formulate two additional models based on different cognitive mapping mechanisms which are similar to those introduced in models \eqref{system-eqn-1}, \eqref{system-eqn-2} and \eqref{system-eqn-3}.
 
 The first extension considers a more realistic scenario where the attractive potential $a(x,t)$ is now dynamic, and will be described as a cognitive map denoted by $q(x,t)$. The map $q (x,t)$ is assumed to grow constantly at rate $b>0$ with respect to resource density, while it decays linearly at a rate $\mu \geq 0$ due to finite memory capacity. The equation describing the evolution of the cognitive map is then
 \eq{
 \frac{\p q}{\p t} &= b v - \mu q ,
 }
 and in \eqref{system-eqn-RC1} we take $a(x,t) = \overline{q}_{g,R}$. Model \eqref{system-eqn-RC1} then becomes
 \eql{
   \begin{cases}
 \frac{\p u}{\p t} = D_1 \D u - \g \grad \cdot \left( u \grad \overline{q}_{g,R} \right)  + \frac{c \b u v}{\a + v} - d u, \\
 \frac{\p v}{\p t} = D_2 \D v + r v (1 - v/K ) - \frac{\b u v}{\a + v}, \\
 \frac{\p q}{\p t} = b v - \mu q .
 \end{cases}
 }{system-eqn-RC1-case2}
 This perspective is comparable to that found in model \eqref{system-eqn-2}: instead of observing marks left on the landscape, consumers are assumed to detect the local resource density and are able to maintain a record of where they have previously found resources.
 
{W}e {also} consider versions of dynamic memory includ{ing} additional mechanisms: $q(x,t)$ is now assumed to grow proportional to the resource density \textit{and} the density of the consumers at rate $b > 0$. This may be more reasonable than the previous model, as one of the implicit assumptions in these memory-based movement models is that foragers are able to share knowledge between individuals. Hence, a location with high resource density is more likely to be remembered if a larger number of foragers perceive it as such. Similar to the previous model, it is assumed that the map decays linearly due to finite memory capacity at rate $\mu \geq 0$, however it is also assumed that the map can decay further at rate $\xi \geq 0$ should the consumer return to an area and find a low resource density. The evolution of $q(x,t)$ is then given by
 \eq{
 \frac{\p q}{\p t} &= b u v - (\mu + \xi u) q.
 }
 Taking again $a(x,t) = \overline{q}_{g,R}$, model \eqref{system-eqn-RC1} becomes
  \eql{
   \begin{cases}
 \frac{\p u}{\p t} = D_1 \D u - \g \grad \cdot \left( u \grad \overline{q}_{g,R} \right)  + \frac{c \b u v}{\a + v} - d u, \\
 \frac{\p v}{\p t} = D_2 \D v + r v (1 - v/K ) - \frac{\b u v}{\a + v}, \\
 \frac{\p q}{\p t} = b uv - (\mu + \xi u) q .
 \end{cases}
 }{system-eqn-RC1-case3}
 This perspective is comparable to that found in model \eqref{system-eqn-3-mod2}, and features similar mechanisms found in model \eqref{system-eqn-RC1-case2}. In this way, it is similar to the memory of direct animal interactions, but the direction of bias is opposite: consumers remember areas they are attracted to, not areas that they seek to avoid.

 \begin{oppr}
In what sense do solutions exist solving the time-dependent problems \eqref{system-eqn-RC1-case2} or \eqref{system-eqn-RC1-case3} subject to periodic boundary conditions with a top-hat detection function in a bounded domain? In what sense do spatially non-constant steady state solutions exist? What qualitative properties do these non-constant steady state solutions hold?
 \end{oppr}

 \begin{oppr}
How do the non-constant steady state solutions differ qualitatively among formulations \eqref{system-eqn-RC1-case1}, \eqref{system-eqn-RC1-case2}, and \eqref{system-eqn-RC1-case3}?
 \end{oppr}
 
 We conclude this section with a brief but important note concerning attracting and repelling quantities in relation to dynamically changing cognitive maps. Since each of the equations introduced to describe the cognitive map should satisfy a positivity lemma (i.e., {they are all linear differential equations} with respect to the cognitive map variable), \textit{a single equation is insufficient to include both attracting and repelling quantities.} Instead, they can only describe relative amounts of attraction \textit{or} repulsion with respect to the variable being remembered, but not simultaneously. We explore this point further in Section \ref{sec:shortlongmemory}.
 
 \subsection{Explicit Memory}\label{sec:explicitmemory}
 
{T}hus far, we have seen perception and some forms of memory, yet these forms are all implicit in that they do not feature an explicit reference to previous experiences. A more recent consideration, strongly motivated by the influence of memory on animal movement, is the inclusion of time delays \textit{in the advective potential}. In this way, foragers make explicit reference to their previous experiences. This {sometimes} complicates the mathematical analysis, which we discuss in Section \ref{sec:biologicalinsights}. {This} increase in complexity is not {entirely} unexpected and these models often yield rather interesting results, both mathematically and ecologically. Nevertheless, in comparison to the models introduced in Sections \ref{sec:staticmemory}-\ref{sec:implicitdynamicmemory}, the overarching theme remains the same: agents have a bias in their movement, and this bias is included within the advection term. The key difference is that the bias is now driven by information from the past. This is perhaps the most explicit inclusion of memory within deterministic population models, as there is a continual, explicit reference to previous experiences or information. Most of the existing models consider the discrete delay case, which features reference to information exactly $\tau>0$ {earlier}, but some also consider the more complicated case of a distributed delay, which includes a reference to \textit{all} previous times under some weighting function. While the latter case is more logical, the former is often valued for its simplicity. {It is worth noting that, unlike the previous sections, many models featuring time delays are purely local, at least in the advection term. For this reason, the well-posedness of these problems is generally easier to obtain by existing techniques. Instead, the focus is on bifurcation analysis and an attempt to more fully describe the potential for spontaneous pattern formation and the qualitative features of non-constant steady states.}
 
 \subsubsection{Discrete Time Delays}\label{sec:discretedelays}

\noindent\textbf{Scalar Equations}

We continue with our motivation from the prototypical diffusion-advection model \eqref{prototype}. Perhaps the most obvious way to include the memory of past occurrences is to consider the case where the advective potential $a (x,t)$ is given by the solution itself, $\tau$ units in the past, i.e., $a(x,t) = u(x,t-\tau)$, which {we denote} by $u_\tau := u(x,t-\tau)$. This is the form derived and investigated in \cite{Shi_JDDE} {in addition to} a birth/death process. The equation describing the evolution of the population $u(x,t)$ is given by
\eql{
\frac{\p u}{\p t} = d_1 \D u - \g \grad \cdot ( u \grad u_\tau ) + f(u),
}{delayprototype}
 where $d_1 > 0$, $\g \in \mathbb{R}$, and the growth term $f(u) \in C^{1} (\mathbb{R}^+)$ is roughly of logistic type{, i.e., $f(0) = f(1) = 0$ and $f(u) < 0$ for $u > 1$}. We again emphasize the relation between the standard diffusion-advection equation \eqref{prototype} and the form found in equation \eqref{delayprototype}: the bias in movement is given explicitly by $u_\tau (x,t)$, where $\tau > 0$ can be thought of as the \textit{averaged} spatial memory period. Intuitively, it {is not reasonable} to assume that organisms make a constant reference to information obtained \textit{exactly} $\tau$ time units ago. Indeed, higher forms of memory are expected to be more complicated than this in reality; {still}, it is a useful starting point to consider the effects of an explicit reference to the past. Similar to model \eqref{scalar-eqn-1-1}, the advection rate $\g$ may have a different sign depending on the situation: $\g<0$ represents a movement away from areas of high population density $\tau$ time units ago, which is a natural phenomenon; on the other hand, $\g>0$ represents a movement towards high population densities $\tau$ time units ago, which may be the case for animals that aggregate for group defence (see the discussion of \cite{Shi_JDDE}).

 \begin{oppr}
In \cite{Shi_JDDE}, a homogeneous Neumann boundary condition is chosen for model \eqref{delayprototype}. How does the analysis performed in \cite{Shi_JDDE} change with respect to changes in boundary data, such as a periodic boundary condition or homogeneous Dirichlet boundary condition?
 \end{oppr}

 \begin{oppr}
How does the inclusion of nonlocal perception change the analysis performed in \cite{Shi_JDDE}? In particular, how does the combination of a single, discrete time delay paired with a top-hat detection function change the dynamical outcomes? By the same reasoning used in Section \ref{sec:implicitdynamicmemory}, it may be easiest to study this in a one-dimensional spatial domain subject to periodic boundary conditions so that no further modification of the detection function near the boundary is necessary. 
 \end{oppr}
 
 This model has since been extended in three ways, with focus given to {modification of} the growth term: first, a nonlocal spatial effect is considered in the growth term; second, a nonlocal temporal effect is considered in the growth term; third, a combination of both of these effects is considered in the growth term. 
 
The first case, {explored in} \cite{Song_JDE}, {considers the equation}
\eql{
\frac{\p u}{\p t} = D_1 \frac{\p ^2 u}{\p x^2} + D_2 \frac{\p}{\p x} \left(u \frac{\p u_\tau}{\p x}\right) + f(u, \overline{u}),
}{delayprototypenon-local}
where $\overline{u} (t) = \as{\O}^{-1} \int_\O u (y,t) dy$ is the average population density over the entire domain. Notice that, in comparison to the averaging with respect to a perceptual kernel, this average may change with respect to time but remains {constant in} space. {This form is} due to a recognition that birth/growth/death rates almost certainly depend on population densities at other spatial locations, not purely on the single point where the organism is located. This particular form of nonlocal (in space) interaction is {inspired} by \cite{Furter1989}, which considers this to be the most straightforward way to include a nonlocal (spatial) interaction effect. {Since the referenced work does not prove the existence of solutions to this problem, we highlight the following open problem.}

\begin{oppr}
Under what conditions does there exist a unique solution solving problem \eqref{delayprototypenon-local} subject to a homogeneous Neumann boundary condition? What about in higher spatial dimensions or other boundary conditions?
\end{oppr}
The second case, considered in \cite{Shi_Nonlinearity}, has a comparable form:
\eql{
\frac{\p u}{\p t} = D_1 \D u + D_2 \grad \cdot \left(u \grad u_\tau \right) + f(u, u_\s),
}{delayprototypenon-localtime}
where $u_\s = u(x,t-\s)$. Biologically speaking, this accounts for a delay in the renewal of resources or the time necessary for animals to reach maturity. Readers should note that while the reaction terms may look similar between models \eqref{delayprototypenon-local} and \eqref{delayprototypenon-localtime}, their interpretation is distinct and may take significantly different forms. 

In the discussion of \cite{Shi_Nonlinearity}, it is noted that the effects of diffusion and time delays are not independent of each other; individuals located at $x$ at a {previous} time may move to a new location at the present time. As a reasonable revision, the third case considers a combination of the previous two scenarios, as done in \cite{An_DCDS}. Their evolution equation takes the form
\eql{
\frac{\p u}{\p t} = \D u + d \grad \cdot \left(u \grad u_\tau \right) + \l u F(u, \overline{u_\s}),
}{delayprototypenon-localtime2}
where $d \in \mathbb{R}$ and $\overline{u_\s} (x,t) = \int_\O K(x,y) u(y,t-\s)dy$ for some reasonably smooth spatial kernel $K(x,y)$ (see \textbf{H2} in \cite{An_DCDS}). The spatial scaling is chosen such that $d$ represents the ratio of the memory-based advection coefficient to the regular diffusion coefficient, and $\l>0$ is a scaled constant. The kernel $K(x,y)$ accounts for the nonlocal intraspecific competition of the species for either resources or space. Readers should note that this kernel is distinct from the perceptual kernels introduced previously, and so we denote it by $K$ instead of $g$. The specific choice in function $F$ depends greatly on the purpose of application, and so interested readers are directed to \cite{An_DCDS, Britton1990, Chen2016} for further details. Lastly, notice that \eqref{delayprototypenon-localtime2} is somewhat a generalization of the simpler form introduced for model \eqref{delayprototypenon-local}, where one recovers \eqref{delayprototypenon-local} by taking $K(x,y) = \as{\O}^{-1}$ in \eqref{delayprototypenon-localtime2}. {In the referenced work, the authors explore the existence of steady state solutions only, which gives the following open problem.}

\begin{oppr}
Under what conditions does a unique solution exist solving the time-dependent problem \eqref{delayprototypenon-localtime2} subject to homogeneous Dirichlet boundary conditions? 
\end{oppr}

\noindent\textbf{Systems of Equations}

{The} recent paper \cite{Song2021Preprint} considers a consumer-resource model, similar to model \eqref{system-eqn-RC1}, {with time-delay in the advective potential}. The resources are assumed to be plants or ``no brainer" animals, so that the prey have no memory-based movement. So far, this is the only model to include an explicit memory mechanism within a system of {PDEs}. {For $\tau>0$, }the {system} in one spatial dimension {is} given as
\eql{
\begin{cases}
\frac{\p u}{\p t} = d_{22} \frac{\p^2 u}{\p x^2} - d_{21} \frac{\p}{\p x}\left( u \frac{\p v_\tau}{\p x} \right) + f(u,v), \\
\frac{\p v}{\p t} = d_{11} \frac{\p ^2v}{\p x^2}+ g(u,v), 
\end{cases}
}{delaymodelsystem}
$v(x,t)$ is the density of the resource, diffusing at rate $d_{11}>0$, and {grows/decays according to} $g(u,v)$. $u(x,t)$ is the density of the consumer, diffusing at rate $d_{22}>0$, and {grows/decay according to} $f(u,v)$\footnote{{To save confusion, readers should note that the role of $u$ and $v$ found here are opposite to that found in the original reference to keep the current work as consistent as possible}}. The consumer $u(x,t)$ then moves up the gradient of the resource $v(x,t)$, $\tau$ units of time ago, at rate $d_{21} \geq 0$. Note that $d_{21}$ is taken to be non-negative since it is assumed that the consumer is attracted to the resource. This {form} is comparable to the resource-consumer model found in \eqref{system-eqn-RC1}, {with the key difference being the cognitive map}. {Existence of solutions is not established, and so the following is an open problem.}

\begin{oppr}
Under what conditions does a unique solution exist solving problem \eqref{delaymodelsystem} subject to homogeneous Neumann boundary data? What about in higher spatial dimensions? 
\end{oppr}

Finally, the competition-diffusion model in \cite{Shi2021a} features memory-based self-diffusion and cross-diffusion. {Similar to previous models}, {for $\tau>0$} the {system} takes the form: 
\eql{
\begin{cases}
\frac{\p u}{\p t} = D_1 \D u + D_{11} \grad \cdot ( u \grad u_\tau) + D_{12} \grad \cdot (u \grad v_\tau) + f(u,v), \\
\frac{\p v}{\p t} = D_2 \D v + D_{21} \grad \cdot ( v \grad u_\tau) + D_{22} \grad \cdot ( v \grad v_\tau ) + g(u,v),
\end{cases}
}{delaymodelcompsystem}
The authors investigate the impact of memory-based self- and cross-diffusion by carrying out a stability and bifurcation analysis. In {this} general form, it includes each form of attraction/repulsion to/from their same group or their competitor, depending on the sign of $D_{ij}$. Due to the complexity of the analysis involved, some simplifications and specific cases are considered for clarity: a Lotka-Volterra competition is investigated, i.e.,
\eq{
f(u,v) = u (1 - u - \a v), \quad g(u,v) = \g v ( 1 - \b u - v),
}
for $\a, \b, \g  > 0$. The local stability of steady states is then explored in relation to the kinetics system (i.e., no diffusion or advection) with a focus on the cases a) $D_{12} = D_{21} = D_{22} = 0$, $D_{11} \neq 0$, or b) $D_{11} = D_{22} = D_{21} = 0$, $D_{12} \neq 0$. Case a) has self-aggregation for species $u$, $\tau$ units in the past, while case b) has attraction or repulsion of species $u$ to/from species $v$, $\tau$ units in the past. When they consider the weak competition case, i.e., $\a \b < 1$, some interesting insights are obtained. First, if $u$ is a timid competitor who moves away from the previous locations of its competitor ($D_{12} > 0$), then the constant coexistence steady state will be destabilized as the rate $D_{12}$ increases. On the other hand, if $u$ is an aggressive competitor who moves towards the previous locations of its competitor, then a Hopf bifurcation occurs as the memory period $\tau$ increases. This indicates that memory-based cross-diffusion alters the {monotone} dynamics of {classical 2-species} competition-diffusion systems {where no such Hopf bifurcation can occur}.\\

\noindent\textbf{Extensions of Discrete Delay Models}

As a logical extension of each of the models above, one may incorporate both memory and perception through a modification of the advective potential $u_\tau (x,t)$. Instead, one may consider $\overline{u_\tau}_{g,R} (x,t)$ for some detection function $g$ and perceptual radius $R$. {We highlight some interesting open problems in this regard now.}

\begin{oppr} 
In what sense do solutions exist to solve problems of form \eqref{delayprototype} or \eqref{delaymodelsystem} when a nonlocal perception through a top-hat detection function is included in addition to an explicit time delay? Similar to problems without time delays, it is likely easier to study these problems under a periodic boundary condition.
\end{oppr}

\begin{oppr} 
In what sense may these problems converge in the limit as $R \to 0^+$? What are the qualitative differences between solutions to local and nonlocal problems featuring time delays? 
\end{oppr}

\subsubsection{Distributed Time Delays}\label{sec:distributeddelays}

As previously suggested, {rather than a discrete delay,} it is more realistic to consider a distribution over all previous times, though it becomes significantly more technical than each of the previous models \cite{Shi_JMB}. {Models with distributed delays following Gamma distributions are equivalent to a system of PDEs with the simplest case equivalent to chemotaxis models (see the Appendix), and thus efforts are perhaps best directed to exploring the properties of these equivalent PDE systems without delay.} Before introducing the model itself, it is instructive to first explore the following \textit{spatiotemporal convolution kernel}:
\eql{
v_u (x,t) = \mathcal{G} * * u(x,t) := \int_{-\i}^t \int_\O G(d_3, x, y, t-s) \mathcal{G} (t - s) u(y,s) dy ds .
}{percepspatim}

{We first describe the motivation for this form}. {I}t is assumed that there is some biological reason for the inclusion of a time delay, which in our case is motivated by an explicit memory-driven movement mechanism and a reference to previous experiences. More precisely, we assume that the cognitive map (in this case, the population density $u(x,t)$) at a time $t$ has contribution from itself at all previous times $s<t$, but not all previous times are equally important. The choice in temporal kernel $\mathcal{G}(\cdot)$ then describes the weighting given to previous times, similar to how the perceptual kernel $g(\cdot)$ describes the weighting given to {different} locations. One may even choose $\mathcal{G} = \d (t-\tau)$, $\tau>0$, if one assumes that the population density exactly $\tau$ time units ago is most important. In order to determine the contribution from all previous times, we then multiply the density $u(x,\cdot)$ at time $s$ by the weight function at this time, which is $\mathcal{G} (t-s)$ since it is $t-s$ time units ago. Then, the integration over space is meant to capture perceptual influences. Different from previous examples, the perception kernel depends on space as well as time, and $G(d_3,x,y,t)$ is chosen to be the fundamental solution to the heat equation with boundary conditions identical to those prescribed in the proposed movement model. {Mathematically, the parameter $d_3$ is the diffusion rate of this fundamental solution, which can theoretically be different than the diffusion rate $d_1$ of the species considered. Biologically, the parameter $d_3$ can be thought of as a \textit{smearing} of the cognitive map as a result of imperfect recollection of spatial information.} Such a choice is due to mathematical convenience at the cost of biological realism.

 In \cite{Shi_JMB}, the waning of memory due to the passage of time is considered through a Gamma distribution, with two specific cases referred to as a \textit{weak} or \textit{strong} kernel, respectively: 
\eql{
\mathcal{G}_w (t; \tau) = \tau^{-1} e^{-t / \tau}, \quad\quad \mathcal{G}_s (t; \tau) = t \tau^{-2} e^{-t / \tau}, \quad \tau > 0.
}{weakstrongkernel}
In {applications}, these two temporal kernels have different interpretations. In {our context}, the weak kernel represents knowledge loss only due to waning memory, {whereas} the strong kernel describes both knowledge gain due to learning and loss due to waning memory. As $t$ increases, the weak kernel $\mathcal{G}_w$ is monotonically decreasing, while the strong kernel $\mathcal{G}_s$ is increasing first and then decreasing. The cognitive map is then $a(x,t) = \g v_u (x,t)$, and evolution of motion with a birth/death process is described by
\eql{
\frac{\p u}{\p t} &= d_1 \D u - \g \grad \cdot ( u \grad v_u ) + f(u).
}{delayfinal}
Here, the advective potential is $v_u (x,t)$ {given by} \eqref{percepspatim} with either the weak or strong kernel defined in \eqref{weakstrongkernel}. Although the mathematical complexity is increased, the primary motivation remains the same: agents in population $u$ have a bias in their movement, and this bias is explicitly driven by an attraction (or repulsion) to (from) all previous locations at rate $\g$, with favour given to experiences closer in space and more recent in time. Model \eqref{delayfinal} can be thought of as the prototypical animal movement model which explicitly includes a distributed memory, in the same way that model \eqref{delayprototype} can be thought of as the prototypical animal movement model which explicitly includes memory through a discrete delay. 

This prototypical model can be {extended in a number of ways}. As model \eqref{delayprototypenon-localtime} was a generalization of \eqref{delayprototype} in the discrete delay case, \cite{Song2021} generalizes \eqref{delayfinal} to include a distributed delay in the maturation process. To this end, {for $i=1,2$,} define
\eql{
v_i (x,t) := \mathcal{G}_i * * u(x,t) := \int_{- \i}^t \int_\O G(d_3, x,y,t-s) \mathcal{G}_{i} (t-s; \tau_i) u(y,s) dy ds .
}{percepspatim2}
In this construction, there are two different delay values $\tau_i$, where $\tau_1$ is related to the delay for inclusion of memory and learning, and $\tau_2$ is related to the delay due to the maturation process. $G(d_3, x,y)$ is a spatial kernel which is again taken to be Green's function for the heat equation subject to homogeneous Neumann boundary conditions. Similar to model \eqref{delayfinal}, the diffusion coefficient for this spatial kernel $G(d_3,x,y)$ is taken to be identical to $d_1$, the diffusion rate for the random movement of the population. The kernels $\mathcal{G}_i(t;\tau)$ are then taken to be either the weak or strong kernel defined previously. This leads to four possible combinations for each kernel type. The evolution equation is then given as
\eql{
\frac{\p u}{\p t} = d_1 \D u + d_2 \grad \cdot ( u \grad v_1) + f(u, v_2),
}{delayprototypenon-localtime-dist}
where $d_1>0$ and $d_2 \in \mathbb{R}$.  {The referenced work explores a bifurcation analysis, but the existence of solutions is not established.}

\begin{oppr}
Under what conditions does a unique solution exist solving problem \eqref{delayprototypenon-localtime-dist} subject to homogeneous Neumann boundary data? In particular, can a result be proven for a general spatiotemporal kernel that includes Green's function as a special case?
\end{oppr}

Finally, the double convolution kernel introduced in \eqref{percepspatim} can be more generally thought of as
\eql{
v_a (x,t) = \mathcal{G} * * a(x,t) := \int_{-\i}^t \int_\O G(d_3, x, y, t-s) \mathcal{G} (t - s) a(y,s) dy ds ,
}{percepspatim-general}
where $a(x,t)$ is any potential related to environmental covariates. For example, {one could replace} $a(x,t) = u(x,t)$ as in the prototypical model \eqref{delayfinal} {with any of the cognitive maps described previously}. {The kernel} \eqref{percepspatim-general} then describes modifications to the cognitive map $a(x,t)$ with respect to both space and time. A general evolution equation which includes modifications in both space and time is given as
\eql{
\frac{\p u}{\p t} &= d_1 \D u + d_2 \grad \cdot ( u \grad v_a ) + f(u) .
}{percepspatim1-general}

Note that the inclusion of distributed temporal delays is a rather recent development as applied to knowledge-based movement models. {The scalar model \eqref{delayfinal} includes the well-known Keller-Segel chemotaxis model as a special case} (see \cite[Section 2]{Shi_JMB}; the lemma statement is found in the Appendix for completeness). It would be an interesting direction of study to consider cases where $a(x,t)$ is a more general cognitive map, similar to those introduced in Section \ref{sec:implicitmemory}, and to compare the results with and without the inclusion of a distributed temporal delay. This could provide insights into the effect of learning on animal movement models, motivating the following group of open {problems}. {These problems are somewhat less well-defined compared to previous open problems since even the construction of such a model would be new.}

\begin{oppr}\label{broadoppr1}
How might a dynamic cognitive map interact with distributed or discrete time delays? That is, how might the effects of a dynamic map as found in Section \ref{sec:implicitmemory} combined with discrete or distributed delays as in the models introduced above influence the space use outcomes? Is it possible to use a standard boot-strapping method (see, e.g., the proof of \cite[Proposition 2.1]{Shi_JMB}) to prove the existence and uniqueness of solutions to problems with both nonlocal perception and discrete or distributed delays in a bounded domain?
\end{oppr}

 \subsection{Short and Long Term Memory}\label{sec:shortlongmemory}
 
In the derivation of the space use coefficients for model \eqref{prototype} (see Section \ref{sec:spacecoeffder}), it is suggested that there should be multiple \textit{layers} or \textit{channels} which together {comprise} the cognitive map. These layers may {include} a number of important environmental covariates {influential to} animal movement, such as food acquisition, territory defence, or mate finding \cite{Fagan2017}. {So} far, only one layer is incorporated into the cognitive map. {Motivating} the inclusion of at least two layers that work together to inform movement is through distinct \textit{short} and \textit{long} term memory components (sometimes referred to as a ``bi-component" mechanism \cite{Bracis2015, riotte2015memory,van2009memory}). From a biological standpoint, {several} species are believed to rely on short and long-term memory for seasonal or long-distance migration \cite{kitchin2002}. Recently, it has been shown in a stochastic setting that both short and long-term memory layers {are} necessary to produce periodic movement in a periodic environment \cite{Lin2021}. However, such a result may be less insightful in a PDE setting as a periodic environment guarantees the existence of a {nontrivial} periodic solution depending on the sign of the principal eigenvalue to a linearized problem \cite[Ch. 2 \& 3]{Hess1991}.

To {explore} the effects of short and long-term memory, one makes the following set of assumptions: short-term memory has larger decay and uptake rates, while long-term memory has smaller decay and uptake rates. {T}his way, long-term memory takes more time to form but decays less rapidly{, while} short-term memory responds to changes quickly and fades easily. If we denote by $m_s$ and $m_l$ the short and long-term memory, respectively, the evolution of the short and long-term memory layers could be described by
\eql{
\begin{cases}
\frac{\p m_s}{\p t} = \a_s a_s (x,t)  - \b_s m_s , \cr
\frac{\p m_l}{\p t} = \a_l a_l (x,t)  - \b_l m_l  ,
\end{cases}
}{shortlongmemory}
with motivation taken from models \eqref{system-eqn-3} and \eqref{system-eqn-2}. The parameters should be chosen so that $\a_l < \a_s$ and $\b_l < \b_s$ to capture the effect of short versus long-term memory. The functions $a_s, a_l$ should, in general, depend on time and space, and will be associated with whichever environmental covariate is being tracked. For example, one could take $a_s (x,t) = a_l (x,t) = {u(x,t)}m(x,t)$, where $m(x,t)$ is a resource density, {giving}
\eq{
a(x,t) = c_1 m_s (x,t) + c_2 m_l (x,t).
}
The evolution equation describing movement with perception and short and long-term memory is given by
\eql{
\frac{\p u}{\p t} &= d \D u - \grad \cdot \left( u \grad \overline{a}_{g,R} \right),
}{shortlongequation}
for an appropriate detection function $g(\cdot)$ and  perceptual radius $R$. The coefficients $c_i$, $i=1,2$, could be either positive or negative depending on application; in \cite{Lin2021} both coefficients are taken positive so that foragers are attracted to both short and long-term memories,{ while }\cite{Bracis2015} suggests that long-term memory could be an attractive force while short-term memory is repulsive, i.e., $c_1 < 0 < c_2$. \cite{Bracis2015} incorporates an attraction to high-resource areas through long-term memory, while being repelled from areas it has recently been through short-term memory in order to allow resources to replenish. This formulation may be useful {to explore} the effects of how the time since last visiting a location influences movement decisions, such as in relation to resource density \cite{Bracis2015,schlagel2014detecting, Ranc2020} or territory surveillance and prey management \cite{Schlagel2017}. Readers should note that in these references, the models are stochastic or statistical. To our knowledge, \eqref{shortlongmemory}-\eqref{shortlongequation} is the first incorporation of short and long-term memory in a diffusion-advection equation setting. This opens an interesting avenue of study, as this combination of effects allows one to include a prioritization of information through an ordering of advection rates. Logically, it makes sense to prioritize avoiding predators over obtained sustenance, for example. Variable rates of advection depending on other external factors (e.g., hunger or other satisfactions measures, see Section \ref{sec:learning}) provide a more complex mechanism which includes the effect of prioritization. {We highlight two broad open problems regarding short and long-term memory. Similar to Open Problem \ref{broadoppr1}, even a reasonable model formulation from this perspective would be new. As such, more precise research questions similar to those found in Section \ref{sec:implicitdynamicmemory} can be formulated once such a model has been constructed.}

\begin{oppr} 
What role does short vs. long-term memory play in differing modelling scenarios, e.g., the well-posedness of problems or influences on steady state profiles? Are both short and long-term memory components necessary and/or sufficient to predict certain space use phenomena, if ever? That is, does the inclusion of both short and long-term memory produce novel effects not found when only, say, short-term memory is included?
\end{oppr}

\begin{oppr}
Can this idea of short and long-term memory be meaningfully connected to the concept of \textit{time since last visit} (TSLV) as explored in \cite{Schlagel2017}? More precisely, does a model featuring both short- and long-term memory components provide better predictive power than models with only one memory compartment? Such insights would provide an interesting connection between discrete-time probabilistic models, as used in \cite{Schlagel2017}, and continuous-time deterministic models of the sort explored here.
\end{oppr}

\begin{oppr} 
How does the prioritization of information affect space-use patterns or persistence/extinction outcomes? For example, is a given model able to predict or explain mechanistically the importance of avoiding predation before seeking food? What does a model predict if the prioritization is reversed? 
\end{oppr}

\subsection{Learning}\label{sec:learning}

The phenomenon of \textit{learning} is more difficult to quantify in the current setting. This is exacerbated by the fact that there are many different forms of learning, ranging from simple habituation (a change in behaviour through repeated exposure to stimulus) all the way to observational learning (learning through mere observation). If we take the psychologists' definition of learning {given in \eqref{learning-psych}}, one may argue that any model featuring a dynamic cognitive map constitutes learning. For example, if a cognitive map $a(x,t)$ is described by an ordinary differential equation
\eq{
\frac{\p a}{\p t} (x,t) = \a (x,t) - \b (x,t) a (x,t),
}
then the growth term $\a(x,t)$ is the {implicit} learning mechanism. Therefore, one can argue that all models presented in Section \ref{sec:implicitdynamicmemory} that feature an additional equation describing the evolution of a cognitive map include learning. On the other hand, models featuring time delays (Section \ref{sec:explicitmemory}) may or may not include a learning process. Models that feature discrete time delays, for example, do not have a learning mechanism. On the other hand, distributed delays with a \textit{strong} kernel include an implicit learning process. This definition of learning may be too broad for further study, {though}, since it does not allow one to distinguish a genuine ``learning process" from other phenomena. This can be problematic in at least two ways. First, some models feature a cognitive map that lies ``outside" the mind of the foragers, as in model \eqref{system-eqn-2} {using} marks on the landscape. One can conceptually distinguish between a learning process and an external process, but they are described in an equivalent way and so {relative influences} cannot be easily distinguished {mathematically}. Second, this broad definition of learning may create difficulties when {comparing with} {empirical} data: {If one tries to compare a model to data to identify whether a learning process has influenced an observed movement pattern, it is impossible to distinguish if the movement pattern is instead a consequence of some simpler mechanism.}

Alternatively, some consider learning to be markedly different than memory when defined as ``modifications to a forager's behaviour through experience/knowledge acquisition" \cite{ThrunPratt1998}, {which is closer to the task-based definition given in \eqref{learning-task}}, {providing} an avenue to determine whether or not a learning process has occurred. From this {perspective}, none of the models introduced thus far feature learning since their movement mechanism \textit{remains the same for all time}. That is, the movement mechanism is an attraction (or repulsion) from the gradient of the cognitive map, however constructed. This means that a population can never learn to change their {meta-level} behaviour in relation to environmental {changes}. This motivates one to consider the effect of variable rates of attraction or diffusion, {that is,} given an advective potential $a(x,t)$ with rate of attraction $\g$, $\g$ should be allowed to change in sign and/or magnitude as new information is obtained. Indeed, it has been empirically shown that certain factors can increase locomotion activity in some species, {inspiring} the concept of \textit{starvation driven diffusion} \cite{Cho2013}, which incorporates a mechanism for the rate of diffusion to \textit{decrease} when an organism is satisfied with its current environment, and will \textit{increase} when the organism is unsatisfied. To this end, a measure of ``satisfaction" must be introduced, and this quantity can be phenomenologically viewed as a learning mechanism. To this end, a \textit{satisfaction measure} was defined in \cite{Cho2013} {as}
\eql{
s = s(x,t) := \frac{\text{food supply}}{\text{food demand}},
}{satisfactionmeasure}
and so it is assumed that the animals have knowledge of the current food supply and demand. If the resource distribution is given as a function $m(x,t)$, one then has
\eql{
s(x,t) = \frac{m(x,t)}{u(x,t)},
}{satisfactionmeasure-1}
where $u$ is the population density of foragers, and so it is assumed that food demand is related directly to the population density. Notice that $s(x,t)$ is similar to the static cognitive map suggested in model \eqref{scalar-eqn-percapitadensity}, but {is now used} to measure satisfaction {instead of acting as a cognitive map}. If $s>1$, the food supply is larger than the food demand, and so motility decrease{s}. If $s<1$, the food supply is smaller than the food demand and motility increase{s}. Thus, the changes to the rate of diffusion should be a composite function $\omega (s)$, where $\omega (\cdot)$ is small for arguments greater than one and large for arguments less than one, {such as}
\eq{
\omega(s) = \begin{cases} d^+, \quad 0 \leq s < 1, \cr d^-, \quad 1 \leq s < \i , \end{cases}
}
for $0 < d^- < d^+ < \i$. In general, it is suggested that a satisfaction measure should satisfy
\eq{
\omega (s) \nearrow d^+ \text{ as } s \searrow 0^+, \quad\quad \omega (s) \searrow d^- \text{ as } s \nearrow \i,
}
where $d^+$, $d^-$ are the maximum and minimum rates of diffusion. 

This alternative view of learning can be incorporated into the advection term as well: instead of modifications to the rate of diffusion, one could introduce modifications to the advection speed, or even the sign of the advection speed. This may introduce mathematical difficulties, however, as one cannot include the advection speed \textit{outside} of the gradient. Instead, the rate must appear \textit{within} the gradient as is found in the derivation of space use coefficients, see {Section} \ref{sec:spacecoeffder}. More precisely, {if the} cognitive map is given by $a(x,t)$ and the advection rate is given by $\g = \g(x,t,u)$, no longer constant and may depend on the population density $u$ and other environmental factors, the correct form with perception would be $- \grad \cdot ( u \grad (\g(x,t,u)  \overline{a}_{g,R}))$, {rather than} $- \g (x,t,u) \grad \cdot( u \grad \overline{a}_{g,R}))$. Furthermore, while $\omega(\cdot)$ must remain positive when modifying rates of diffusion (in order to maintain parabolicity of the equation), $\g(x,t,u)$ can change sign when modifying advection speeds. In such a case, a sign change would indicate a change from attraction to repulsion (or vice versa), which may be a stronger indication of learning as it indicates a change in \textit{kind} rather than a change in \textit{amount}. Such {forms} may also allow one to overcome an issue discussed in Section \ref{sec:implicitdynamicmemory}: a single equation could, in principle, be sufficient to describe a cognitive map if the rate of advection is allowed to change sign. This point is raised briefly in the discussion of \cite{Shi_JDDE}, where it is suggested that the sign of the advection speed may change from {attraction} to {repulsion} depending on environmental conditions.

To demonstrate this point, we construct a simple model combining the effects of starvation-driven advection and a den site. First, we assume that there is some constant rate of attraction to a central den site located at $x_0\in\O$. We then assume that foragers are attracted to the local resource density, but the rate at which they move up the gradient of the resource profile will now depend on the satisfaction measure $s(x,t)$ defined in line \eqref{satisfactionmeasure-1}. We then define a modified advection rate $\tilde{\omega}(s)$ to be
\eql{
\tilde{\omega}(s) = \begin{cases} \g^+, \quad 0 \leq s < 1, \cr 0, \quad 1 \leq s < \i , \end{cases}
}{SDA}
so that the advection rate depends on the starvation level. Should the forager be hungry, the advection rate ``turns on" at constant rate $\g^+ > 0$; if the forager is not hungry, the advection rate ``turns off" and the default mode tending toward the den site will dominate. With these effects included, the one-dimensional model takes the form
\eql{
\frac{\p u}{\p t} &= \frac{\p}{\p x} \left( d \frac{\p u}{\p x} - u \frac{\p}{\p x}\left[ \tilde{\omega}(s) m {+} \g \norm{x-x_0} \right] \right),
}{SDA-densite-1}
where $0 < \g < \g^+$ so that the movement towards resources becomes a priority when {hunger is high}.

\begin{oppr}
In what sense do solutions exist to solve problem \eqref{SDA-densite-1}? Due to the possibility of the equation remaining linear (in the variable $u$), this should not be as difficult a task as other well-posedness problems discussed thus far in principle. On the other hand, the form described here features irregular coefficients at higher order, which introduces other difficulties.
\end{oppr}

\begin{oppr}
How can one meaningfully incorporate nonlocal perception in addition to the mechanism of learning outlined in model \eqref{SDA-densite-1} and the preceding discussion? What differences in solution behaviour does one observe with only one mechanism (e.g., learning or perception) vs. both mechanisms simultaneously? 
\end{oppr}

Motivated by the formulations above, one may incorporate the idea of a satisfaction measure with those introduced in \cite{Ranta1999} so that foragers may update their advection rates based on (possibly incomplete) knowledge of the local resource density and an \textit{expected} resource density. Recall the form of the cognitive map $m(x,t)/\overline{m}(t)$ introduced in model \eqref{scalar-eqn-averagedensity}, where $\overline{m}(t)$ denotes the average resource density at time $t$. The advection rate can then be taken as a composite function with the new satisfaction measure $s(x,t) = m(x,t) / \overline{m}(t)$, where $s>1$ in above average resource areas and $s<1$ in below average resource areas. Depending on the expected behaviour of the forager, $\omega (s)$ can be formulated so that the advection speed either increases or decreases, or even changes sign, depending on this satisfaction ratio. This begins to address a key difference between spatial memory and attribute memory: the formulations above provide at least one direction to distinguish between the \textit{quality} of movement decisions and behaviour, rather than the mere locations that are desirable or undesirable.

\begin{oppr} 
For the new model \eqref{SDA-densite-1}, can the starvation-driven switch between two cognitive mechanisms, foraging and home attraction, generate nontrivial spatiotemporal dynamics and patterns? Biologically, would such realistic switches benefit or impair the species? What other satisfaction measures can be introduced to act as a behaviour modification mechanism? Due to the potential complexity of the model, even a detailed exploration of simulated solutions for this problem may provide meaningful insights into potential space-use outcomes.
\end{oppr}

\section{{Sample Model Derivations}}\label{sec:derivations}

\subsection{Derivation of Fokker-Planck Equation}\label{sec:FPDer}

In this section, we derive a continuous-time, continuous-space model via the master equation
\eql{
u(x, t+\tau) = \int_\O f(x,y,t ; \tau) u (y, t) dy
}{der1.1}
in an unbounded domain $\O = \mathbb{R}^2$. Higher dimensional derivations can be considered in a similar fashion. The master equation above keeps track of the density function for the location of animals over time and space via a conservation law. The function $f$ is a probability density function which describes the movement from a point $y$ to $x$ over the time interval $[t, t+\tau)$ (see Section \ref{sec:spacecoeffder} for a concrete example). 
 
First, consider the case where an individual is released with its location given by the probability density function $u_0 (x)$ at time $t = 0$. {For example}, one {may} consider the case when an individual is released at a point $x_0 \in \O$ {so that} $u_0 (x) = \d (x - x_0)$. If we do not know the exact location of release, we use a general probability density function (PDF) $u_0 (x)$ which integrates to $1$. We {then} use a Taylor series approach while taking the limit as $\tau \to 0^+$, see \cite{Bharucha-Reid:1960:ETMP}. Starting from \eqref{der1.1}, we define the new vector $z = x - y$. Then we write $f(x,y,t; \tau) =: f_z (z, y,t; \tau)$. Through this change of variables we obtain from \eqref{der1.1}
\eql{
u(x, t + \tau) = \int_\O f_z ( z, x - z, t; \tau) u (x - z, t) dz.
}{a1.2}
Expanding the right hand side in space about $0$, we find
\eql{
u(x, t+\tau) &= \int_\O   f_z (z,x,t;\tau) u(x,t) dz \nn \\
& - \int_\O \left( z_1 \frac{\p}{\p x_1} ( u(x,t) f_z (z,x,t;\tau) ) + z_2 \frac{\p}{\p x_2} ( u(x,t) f_z (z,x,t;\tau) ) \right)dz \nn \\
& + \int_\O \left( \frac{z_1^2}{2} \frac{\p^2}{\p x_1 ^2} (u(x,t) f_z (z,x,t;\tau)) + \frac{z_2 ^2}{2} \frac{\p^2}{\p x_2 ^2} ( u(x,t) f_z (z,x,t;\tau)) \right) dz \nn \\
& + \int_\O \left( z_1 z_2 \frac{\p^2}{\p x_1 \p x_2} ( u(x,t) f_z (z,x,t;\tau) ) \right) dz + \mathcal{O} (z^3),
}{a1.3}
where we have assumed that the mixed partial derivatives of $f_z u$ agree. Using the fact that $\int_\O f_z (z,x,t;\tau) u(x,t)dz = u(x,t)$, we may move this term to the left hand side of \eqref{a1.3} and divide by $\tau>0$ to obtain 
\eql{
\frac{u(x,t+\tau) - u(x,t)}{\tau} &= - \frac{1}{\tau}  \int_\O \left( z_1 \frac{\p}{\p x_1} ( u(x,t) f_z (z,x,t;\tau) ) + z_2 \frac{\p}{\p x_2} ( u(x,t) f_z (z,x,t;\tau) ) \right)dz \nn \\
&+ \frac{1}{\tau}  \int_\O \left( \frac{z_1^2}{2} \frac{\p^2}{\p x_1 ^2} (u(x,t) f_z (z,x,t;\tau)) + \frac{z_2 ^2}{2} \frac{\p^2}{\p x_2 ^2} ( u(x,t) f_z (z,x,t;\tau)) \right) dz \nn \\
&+ \frac{1}{\tau}  \int_\O \left( z_1 z_2 \frac{\p^2}{\p x_1 \p x_2} ( u(x,t) f_z (z,x,t;\tau) ) \right) dz \nn \\
&+ \mathcal{O} ( \tau^{-1} z^3 ) .
}{a1.4}
Taking the limit as $\tau \to 0^+$, we find that
\eql{
\frac{\p u}{\p t} &= - \grad \cdot ( c(x,t) u) + \sum_{i, j = 1}^2 \frac{\p^2}{\p x_i \p x_j} ( d_{ij}(x,t) u) ,
}{a1.4b}
where
\eql{
c(x,t) &:= \lim_{\tau \to 0^+} \frac{1}{\tau} \int_\O z f_z (z,x,t;\tau) dz , 
}{a1.5a}
\eql{
d_{ij} (x,t) &:= \lim_{\tau \to 0^+} \frac{1}{2 \tau} \int_\O z_i z_j f_z (z, x,t; \tau) dz .
}{a1.5b}
Note that we tacitly assume that all higher order terms vanish as $\tau \to 0^+$. Interested readers are directed to \cite{Bharucha-Reid:1960:ETMP} and \cite[Ch. 2.3]{Marley2020} for comparable derivations and further discussion.

\subsection{Derivation of Space Use Coefficients from Utilization Distribution Function}\label{sec:spacecoeffder}

{Next, w}e derive a general form of the space use coefficients $c, d$ from the utilization distribution function
\eql{
f(x,y,t; \tau) = \frac{K ( x-y; \tau) w ( A (x,t))}{\int_\O K( z - y; \tau) w (A(z,t)) dz} ,
}{a1.5.1}
with further considerations found in \cite{Marley2020}; {see also \cite{Potts2022}.} The kernel $K(x-y;\tau)$ is the spatial dispersal found in the absence of external factors (e.g., in a spatially uniform environment), where we assume $K$ is symmetric for simplicity. {The function $A(x,t)$ is a description of relevant environmental covariates, which in our context will seek to include cognitive mechanisms. In general, $A(x,t) = \sum_{j} \b_j a_j (x,t)$ for some list of covariates $a_i (x,t)$ and their associated selection coefficients $\b_j$.} The weighting function $w( \cdot )$ then describes how {such} factors influence the movement behaviour towards the point $x$ at time $t$. This is perhaps the most important term {here}, as this is where we are able to explicitly include cognitive factors from first principles.

Under the same change of variables as in {Section} \ref{sec:FPDer}, the utilization distribution function is given by
\eql{
f_z (z,x,t;\tau) = \frac{K(z;\tau) w(A(z+x,t))}{\int_\O K(\tilde{z};\tau) w(A(\tilde{z}+x,t)) d\tilde{z}},
}{a1.5.3}
and so if one performs another Taylor expansion of $w$ appearing in the numerator of \eqref{a1.5.3}, we {find}
\eql{
K(z;\tau) w(A(z+x,t)) &\sim w(A(x,t)) K(z; \tau) + z \grad w(A(x,t)) +  \mathcal{O} (z^2).
}{a1.5.2}
Hence, under the assumption that the kernel $K$ is symmetric, it is readily seen that
\eql{
\int_\O z_1 z_2 f_z (z, x; \tau) dz &\sim w(A(x,t)) \int_\O z_1 z_2 K(z; \tau) dz + \mathcal{O} (z^3) = 0 + \mathcal{O} (z^3) .
}{a1.6}
Hence, referring to \eqref{a1.5b} one finds that in fact $d_{ij} \equiv 0$ for $i \neq j$. When $i=j$, we further expand the denominator of \eqref{a1.5.3} to obtain
\eql{
\int_\O K(\tilde{z};\tau) w(A(\tilde{z}+x,t)) d\tilde{z} &\sim w(A(x,t)) + \frac{\D w(A(x,t))}{2!} M_2 (\tau) + \ldots,
}{a1.7}
where
\eql{
M_p (\tau) = \int_\O \as{y}^p K(y;\tau) dy
}{a1.8}
is the $p$th moment of the dispersal kernel.

Notice that the first order term in \eqref{a1.7} vanishes in the same way as $d_{ij}$ for $i \neq j$ since $K$ is symmetric, i.e., $M_1 (\tau) = 0$. Similarly, the mixed terms also vanish. Putting expansions \eqref{a1.5.2} and \eqref{a1.7} into \eqref{a1.5b}, we find in the limit that
\eql{
d_{ii} (x,t) &= \lim_{\tau \to 0^+} \frac{1}{2 \tau} \frac{w(A(x,t)) \frac{M_2 (\tau)}{2} + \ldots}{ w(A(x,t)) + \D w(A(x,t)) \frac{M_2 (\tau)}{2} + \ldots} \nn \\
&= \lim_{\tau \to 0^+} \frac{M_2 (\tau)}{4 \tau} + \mathcal{O} (\tau) \nn \\
&= d \ (= constant).
}{a1.9}

Using a similar procedure, one may insert expansions \eqref{a1.5.2} and \eqref{a1.7} into \eqref{a1.5a} to find
\eql{
c(x,t) &= \lim_{\tau \to 0^+} \frac{1}{\tau} \frac{\grad w(A(x,t)) \frac{M_2 (\tau)}{2} + \ldots}{w(A(x,t)) + \D w(A(x,t)) \frac{M_2 (\tau)}{2} + \ldots } \nn \\
&= \lim_{\tau \to 0^+} \left( \frac{M_2 (\tau)}{2 \tau} \frac{\grad w(A(x,t))}{w (A(x,t))}  \right) + \mathcal{O} (\tau) \nn \\
&= 2 d \grad \ln ( w(A(x,t)) )  .
}{a1.10}

Finally, we discuss the choice of weighting function $w$ as it applies to knowledge-based animal movement. Indeed, the final form of the equation describing animal movement is, roughly, a diffusion equation with bias in movement given by advection up the gradient of the log of the weighting function $w$. We now consider the end result under the assumption of \textit{exponential} weighting of covariates, which is the most common choice of selection function appearing in the literature \cite{Fieberg2021,Potts2020}{, including both }environmental {and} cognitive factors influencing movement behaviours. In such a case,
\eq{
w({A(x,t)}) \propto \exp\left(\sum_i \beta_i a_i(x,t) \right) ,
}
where $\beta_i>0$ {($<0$)} indicates attraction {(repulsion)} towards {(away from)} the {environmental covariate} $a_i (x,t)$. Thus, from the derivation of $c(x,t)$ appearing in line \eqref{a1.10}, we see that
\eql{
c(x,t) &= 2 d \grad \ln (w) = 2 d \sum_{i} \b_i \grad a_i (x,t) .
}{a1.11}
Thus, the final description of motion in two spatial dimensions is given by
\eql{
\frac{\p u }{\p t} (x,t) &= d \D u (x,t) - 2d \grad \cdot \left( u (x,t) \left( \sum_{i} \b_i \grad a_i (x,t) \right) \right).
}{FinalFPEqn}

\section{Biological Insights Through Mathematical {Exploration}}\label{sec:biologicalinsights}

In this section, we discuss some of the biological insights we have gained, and continue to gain, through the analysis of the models proposed so far. First we bring attention to some useful ``Rules of Thumb" that {should be considered once} a model has been formulated. {Then, we} discuss {possible} ``Measures of Success" in animal movement models {and} their importance to the {field}. In some cases, we are able to discuss existing insights gained through such measures. In other cases, we propose modifications or suggestions for further study. We conclude with a discussion of {other common avenues of analysis, including pattern formation, travelling wave solutions, critical domain sizes, and useful numerical techniques.} 

\subsection{Rules of Thumb}\label{sec:RoT}

{We first} introduce some rules of thumb that should be considered {for virtually any} model. We divide these rules into ``prerequisites" and ``considerations". Prerequisites are properties that the model should necessarily possess {to remain logically valid.} Considerations, on the other hand, are helpful suggestions that may be {more or less relevant depending on the situation.}

\subsubsection{Prerequisites}

We first discuss {an essential property that} any reasonable model should possess: \textit{existence {of a solution}}. {Without this, any analytical or numerical insights gained may be misleading at best, or entirely incorrect at worst.}

\noindent\textbf{Existence}: Existence simply refers to the existence of a solution to a given problem in some suitable context. It is well known that there are problems for which solutions exist, problems for which no solutions exist, and problems for which we do not know whether a solution exists or not. Furthermore, when the solution depends on a temporal variable, it may exist for all time $t \in (0,\i)$, or it may only exist on some finite interval $(0,T)$ with solution blowup as $t \to T^-$ \cite[Ch. 1.6]{cantrell2004spatial}. 

While it is not inappropriate to investigate ecological models without knowing the existence of solutions a priori (consider the famous Navier-Stokes equations), it is often possible to show the existence of solutions under relatively weak assumptions. Depending on the context, one may obtain \textit{classical} solutions \cite[Ch. 8]{Wu2006Elliptic} (sufficiently differentiable in both space and time), but other contexts may require a more general form of solutions, such as \textit{weak} or \textit{strong} solutions \cite[Ch. 2, 3 \& Ch. 9]{Wu2006Elliptic} or \textit{mild} solutions \cite{Pazy1983,Amann1995}. These are solutions that may not be differentiable in the classical sense, but can be {made} well-defined through notions of weak differentiability or semi-group theory. The technical details are beyond the scope of this review; {we appeal to the references provided above for further details}. It is also worth noting that many of the models proposed here do not fall within the standard theory. This leaves many open problems regarding the existence of solutions. {We have made an effort to highlight explicitly some of these open problems in this regard.}


\subsubsection{Considerations}\label{sec:considerations}

Next, we briefly highlight some of the features anyone should consider when proposing a new model. These properties are by no means necessary for a model to be valid, but they offer insights into the question \textit{under what conditions or assumptions is the model most valid?}

\noindent\textbf{Uniqueness}: The property of uniqueness ensures that there is exactly one non-trivial solution which solves the problem. For time-dependent problems, this implies that there exists exactly one non-trivial solution corresponding to a given initial condition. In terms of ecological (and other) modelling efforts, it is a desirable property to hold to draw concrete conclusions from the model. Indeed, should a model have two or more solutions corresponding to the same set of initial data, it is no longer possible to determine which outcome could be realized in the {natural world}. Furthermore, should a proposed model suffer from a lack of uniqueness, this may indicate that an important biological mechanism has been neglected or ignored, which may motivate alternative model constructions. {This gives a broad open problem to virtually all models introduced in this manuscript.}

\begin{oppr}
The following question can be proposed for any previously introduced model with existence as an open problem: Do there exist necessary and/or sufficient conditions that guarantee the uniqueness of the solution obtained? The results in \cite{Giunta2021,Jungel2022}, for example, prove a uniqueness result, each in a slightly different context.
\end{oppr}

\noindent\textbf{Continuity {with respect to} initial data}: This property roughly says that small changes in the initial data result in comparably small changes in the solution. This is important for modelling purposes particularly when efforts are made towards fitting {empirical} data: small errors in data measurements may produce unreasonable results if this continuity {property} does not hold. {On the other hand, there are examples where small changes \textit{do} result in large changes after a long time period (such as climate and weather models); hence this property is {dependent on context.}}

Uniqueness and continuity {with respect to} initial data, together with the existence of a solution, yield the well-posedness of a problem in the sense of Hadamard \cite[Ch. 1.3.1]{evans1998partial}, which gives a reasonable expectation that the solution can be solved for using standard numerical methods. On the other hand, if the problem is ill-posed in some way, alternative numerical methods may be necessary in order to fully justify the result(s) one obtains. {For example, there are problems for which no solutions exist while a numerical solver will produce a ``solution" that is, essentially, a numerical artifact (see Section \ref{sec:travwave} for a precise reference).} {Another} example is in model \eqref{system-eqn-3} for 2 interacting populations: {through} a linearization about constant steady states, the problem becomes ill-posed as the perceptual radius $R$ tends to zero. {Solutions can be found} at arbitrarily high wavenumbers as a consequence \cite{Potts2016}. {This suggests that 1. the solution profile at steady-state fails to be unique, and 2. the problem almost certainly does \textit{not} enjoy continuous dependence with respect to initial data.}

\noindent\textbf{Initial Conditions:} We briefly {reflect} on the nature of the initial {data} chosen. First, {we note} that the {choice of} initial data can be {chosen almost} arbitrar{ily} in theory {(continuity of the initial data is often sufficient to ensure well-posedness of the problem, for example)}; the importance, perhaps, is on determining the \textit{impact} that differing initial data may have on the {dynamics predicted by the model}. In many classical models, initial conditions have {little} impact on {the overall dynamics} (e.g., the standard Fisher-KPP equation \cite[Ch. 4.1]{Perthame2015}); on the other hand, some models have significantly different outcomes depending on the initial data chosen (e.g., the Allen-Cahn equation \cite[Ch. 4.1]{Perthame2015}). It is not clear how the choice in initial data affects space use outcomes in a general sense, especially in relation to the models introduced here. However, there is strong numerical evidence to suggest that patterns may only form if the initial condition is sufficiently close to a segregation pattern, see \cite{Potts2016}. {Pattern formation should therefore not be considered ``spontaneous" in the {usual} sense.} Instead, such outcomes may be {better} thought of as ``pattern \textit{stimulating}". {O}ne exception to {an} arbitrary choice in initial data is in relation to {empirical data}: should a model be constructed {to fit} {empirical} data, the initial data {should} match that of the {emperical }{data} itself {in absence of further justification.}

\noindent\textbf{Boundary Conditions:} Modellers spend a significant amount of time describing what happens within the domain (e.g., the habitat), but we must also specify what occurs when an agent reaches the boundary of the domain (e.g., the edge of the habitat in a bounded domain). Often this is presented from a mathematical perspective with {three} main types (Dirichlet, Neumann, Robin) introduced. {For application purposes, we also discuss zero-flux and periodic boundary data.} {Here, w}e will approach this first from the biological perspective, and then make some brief comments on the connection to prospects of mathematical analysis.

To motivate the reader, we consider some of the models introduced thus far. Many feature no birth {or} death processes (see Sections \ref{sec:perception}-\ref{sec:implicitdynamicmemory}, for example), and so it is reasonable to assume that the total population remains fixed for all time{;} the models proposed describe \textit{animal movement only}. {Hence}, the most appropriate boundary condition is referred to as the \textit{zero-flux} or \textit{reflecting} boundary condition{, which conserves the total population through construction.} To determine the appropriate zero-flux boundary condition, one simply integrates the equation over the domain and applies the divergence theorem to show that $\frac{d}{dt} \int_\O u(x,t)dx = 0$, which naturally implies {a conserved population}, and is identical to the initial population {size}. As a simple but instructive example, consider \eqref{prototype} in $\O = (0,L)$. Integrating over $(0,L)$ and applying the divergence theorem yields
\eq{
0 = \frac{d}{dt} \int_0 ^L u(x,t) dx &= \int_0 ^L \frac{\p }{\p x} \left(d  \frac{\p u}{\p x} - \frac{\p a}{\p x} u \right) dx = d \frac{\p u}{\p x}   - \frac{\p a}{\p x} u \biggr\vert_0 ^L .
}
Hence, the total population is conserved if we prescribe at the endpoints
\eql{
d \frac{\p u}{\p x}   (0,t)  - \frac{\p a}{\p x} (0,t) u(0,t) = d \frac{\p u}{\p x} (L,t)  - \frac{\p a}{\p x} (L,t) u(L,t) = 0.
}{zerofluxboundary}

We {connect this idea to} another common boundary condition, referred to as a homogeneous Neumann boundary condition, which prescribes the outer normal derivative of the solution at the boundary:
\eq{
\frac{\p u}{\p \mathbf{n}} (x,t) = 0, \quad \text{ on $\p \O$.}
}
Here, $\frac{\p}{\p \mathbf{n}}$ denotes the outer unit normal vector to the boundary of the domain $\p \O$. Sometimes, a homogeneous Neumann boundary condition is also referred to as a zero-flux boundary condition. This raises an interesting point that may contribute to confusion. In some cases, a homogeneous Neumann boundary condition is sufficient to ensure a conserved population{, and so} is equivalent to a zero-flux condition. This would be the case in any model where the potential $a(x,t)$ depends {linearly} on the solution $u(x,t)$ itself {(e.g., {the local version of} model \eqref{scalar-eqn-1-1})}. In this way, $\frac{\p a}{\p \mathbf{n}} = \frac{\p u}{\p \mathbf{n}} = 0$ along $\p \O$, and so \eqref{zerofluxboundary} is automatically satisfied. On the other hand, a homogeneous Neumann boundary {condition} may not conserve the population if $a(x,t)$ is given a priori {and does not itself satisfy a homogeneous Neumann boundary condition}. {In such cases, a Neumann condition is not equivalent to a zero-flux condition.}

The {third} boundary condition we discuss is the homogeneous Dirichlet boundary condition:
\eq{
u(x,t) = 0, \quad \text{ on $\p \O$.}
}
This is sometimes referred to as a hostile boundary condition. As this name suggests, this condition assumes that the boundary is completely lethal and any agent that reaches the edge is removed immediately and never returns (e.g., the animal dies). Such a condition is most appropriate for species on an island, surrounded by a cliff edge, or in a particular necessary ecological niche, for example. In general, a homogeneous Dirichlet condition will not conserve the population since a hostile boundary is a mechanism by which there is a loss of the population without an explicit birth or death process included in the model equation.

{The fourth condition, easily viewed as a generalization of the Neumann and Dirichlet condition,} is {often} referred to as a \textit{Robin} or \textit{mixed} boundary condition. This {is usually} written as
\eql{
 \a (x) \frac{\p u}{\p \mathbf{n}} + \b (x) u = 0, \quad \text{ on $\p \O$,}
}{robinboundary}
but can be generalized to include cases where the coefficients are time dependent \cite[Ch. 2]{pao1992nonlinear}. {When $\a=0 < \b$, it is a Dirichlet condition; if $\b = 0 < \a$, it is a Neumann condition.} If $\a, \b > 0$, this boundary condition can be interpreted as a partial loss of the population at the boundary. Choosing $\a(x,t) = d$, $\b(x,t) = - \frac{\p a}{\p x}$, the Robin condition \eqref{robinboundary} is equivalent to the zero-flux condition \eqref{zerofluxboundary}. 

{The final boundary condition we discuss in detail is the periodic boundary condition. As the name suggests, this condition produces a level of continuity from one boundary portion to another, giving an implicitly defined ``periodic extension" of a solution from a finite to an infinite domain. In one spatial dimension $(0,L)$ this takes the form
$$
u(0,t) = u(L,t), \quad u_x (0,t) = u_x(L,t).
$$
In higher dimensions, this is equivalent to studying a given problem on a torus. Interestingly, periodic boundary conditions often preserve the total population in absence of birth/death processes for many models introduced here. An unfortunate reality, perhaps, is that while a periodic boundary condition may not be the most biologically reasonable choice, it is much easier to study analytically. The reason is easy enough to understand conceptually: the periodic boundary condition case is the \textit{no boundary} case. For nonlocal problems, this is of particular interest {as} it allows one to avoid addressing what happens to the perceptual kernel near the boundary of the domain. Instead, the kernel is able to ``spill over" and remains a well-defined quantity throughout the entire domain. To observe this visually, consider the third panel of Figure \ref{StaticCogMap}: what happens when the centre of the ball is chosen to be, say, $(1,1)$? For periodic conditions, the portion that spills over in the bottom left corner appears in the top right corner. This is in contrast to the strategy found in \cite{Potts2016}, for example, where the perceptual kernel is explicitly defined as
\eql{
\overline{a}(x,t) := \begin{cases} \frac{1}{R + x} \int_{-x}^R a(x+z,t) dz \quad &\text{ if } 0 < x \leq R, \cr
\frac{1}{2R} \int_{-R} ^{R} a(x+z,t) dz, \quad &\text{ if } R < x < L-R, \cr
\frac{1}{R + L - x} \int_{-R} ^{L-x} a(x+z,t)dz, \quad &\text{ if } L - R \leq x < L,\end{cases}
}{cutoffkernel}
where $a(x,t)$ is {the} cognitive map and $0 < R < L$ is the perceptual radius. Instead of spilling outside of the domain, the kernel is now defined in such a way that this issue is removed by explicitly cutting off the kernel near the boundary. This definition is perhaps the most biologically reasonable and is compatible with a no-flux boundary condition; {yet}, models with such a boundary condition and nonlocal component(s) are difficult to study, both analytically and through simulation (see Section \ref{sec:numericalanalysis} for relevant numerical techniques and challenges). As such, existing tools need to be modified, or new tools need to be developed, to fully explore such problems. Even questions of the existence of solutions can be highly non-trivial based on seemingly minor changes to prescribed boundary data. Because of this trade-off, most open questions proposed so far focus on cases subject to a periodic boundary condition in a bounded domain. These cases allow researchers to {most easily} focus on the influence of varying cognitive influences. Hence, we introduce the following open problem that essentially multiplies the number of open problems by the number of potential boundary conditions different from the periodic case.
}

\begin{oppr}
Under what conditions do solutions exist to problems introduced in Section \ref{sec:implicitdynamicmemory} when subject to a homogeneous Neumann or Dirichlet boundary condition and the perceptual kernel defined in the style of \eqref{cutoffkernel}? How do these results compare to the periodic boundary case? 
\end{oppr}

{Moreover, we propose the following general open problem to expand upon the existing definitions of perceptual kernels (one appropriate for a periodic boundary condition or an unbounded domain, another which explicitly cuts off the kernel near the boundary). There have been strides in this direction for nonlocal models different from those proposed here, such as cell adhesion models (see \cite{Hillen2021} and \cite{Hillen2009} and the references therein). {While non-trivial} in general, a clever shift in perspective may alleviate technical difficulties in ways that are unexpected.}

\begin{oppr}
What are some possible alternative definitions or constructions of perceptual kernels that are both biologically reasonable while lending themselves to rigorous mathematical analysis? In particular, how can one reasonably deal with such nonlocal components near the boundary of the domain?
\end{oppr}

{While classical forms of boundary data feature many existing tools well-developed for rigorous analysis,} a variety of {novel} boundary conditions should {also} be considered for biological purposes and application. {One example {not explored here} are} moving boundary conditions, see e.g., \cite{feng2021fisherkpp}, where a non-Stefan free boundary condition was proposed via a trade-off between shorter and longer spatial scales. 

In {some} sense, there are no correct or incorrect choices {in boundary data} should sufficient justification be provided. {As one may observe from the preceding discussion, though, some boundary data may lend themselves to a more rigorous analysis than others; some boundary data may be more biologically reasonable than others. Care should be taken on {either} front.}


\noindent\textbf{Domains:} Related to the choice of boundary conditions is the domain in which the model is studied. In most models introduced here, a one-dimensional spatial domain is chosen for ease of analysis {and model conception} (exceptions include some of the time-delay models introduced in Section \ref{sec:explicitmemory} where increasing the spatial dimension essentially leaves the analysis unchanged). {For such one-dimensional problems,} this can be in a bounded or an unbounded interval. {As discussed above, appropriately defining a nonlocal component near a fixed boundary point is not a simple task. However, studying a nonlocal problem in an unbounded domain removes the necessity to modify the definition of a perceptual kernel. In this sense, an unbounded interval is comparable to studying the problem in a finite interval with periodic boundary conditions. Furthermore, an unbounded interval lends itself to perceptual kernels that do not have a finite domain of definition. This is the case for the exponential and Gaussian kernels defined in \eqref{detectionkernels}.}

\noindent\textbf{Timescales:} A significant portion of movement models described here do not include a birth/death process. In practice, this may be reasonable if one assumes that animal movement occurs at a timescale much faster than the birth/death processes \cite{Potts2019}. However, this may result in {other} difficulties. For example, a common technique is to take advantage of the so-called \textit{quasi-steady state approximation} in the equation for the cognitive map. That is, if a map $a(x,t)$ is described by a dynamic equation, one sets $\frac{\p a}{\p t} = 0$ and can then {(in principle)} solve for $a(x,t)$ explicitly {(in terms of the population density variable $u$)}. This reduces the number of equations, but it may not be appropriate to make such a simplification if the growth of a cognitive map occurs at a timescale that is the same as that of the actual movement. This may {not} be the case for model \eqref{system-eqn-2}, where the ``cognitive map" is given by marks on the landscape, which necessarily occur at the same timescale as the animal movement. {Subsequent} biological insights may not be accurate {as a consequence}. On the other hand, since perception occurs very quickly and influences movement for a longer period of time, it may be reasonable to use a quasi-steady state approximation for models which describe the cognitive map itself, such as model \eqref{system-eqn-3}. {An alternative} is to include a birth/death process in the model at the start so that the timescale over which these models are valid becomes less of an issue. {Even still, one should take care in including birth/death processes in a model when drawing connections to data collected over a timescale that does not capture significant changes in the population.}

We wish to emphasize that none of these rules are always wrong or right. Instead, we hope to {encourage care in using} some of the assumptions commonly made, {recognizing} when they are most valid, {or} when they may need to be reconsidered{/reformulated}. 

\subsection{Measures of Success in Movement Models}\label{sec:measuresofsuccess}

{In the development of movement models with cognitive mechanisms, it is of great interest to be able to evaluate and compare the outcomes of different formulations. Indeed, it is possible to describe, either quantitatively or qualitatively, the asymptotic dynamics or the possible forms of steady states for a given problem; however, such insights occasionally remain a mere mathematical description rather than a more explicit biological one.} Having a {clearer}, unified framework to measure the {utility} of certain cognitive mechanisms {in a given \textit{biological} context} will be invaluable to the assessment of {model validation and prediction(s)}. In this sense, a \textit{unified} framework may be interpreted more accurately as \textit{unified within groups of comparable models}. {To clarify this point, we briefly describe some of the existing measures used in animal movement models, which suggests in what sense models may be ``comparable".}

{One of the most common measures of success used in population dynamics {generally} is the concept of \textit{fitness}. Classically, fitness corresponds directly with the ability of an organism or population to successfully reproduce. For the majority of cases considered here, we study a population density or a population's probability density function, and so we must study the fitness of the population. In Section \ref{sec:popmeasure}, we discuss measures of success with respect to this sense of fitness. In practice, it is mathematically and conceptually most simple to define a measure of fitness for models that include explicit population dynamics. It is in this sense that such a group of models can be referred to as comparable.}

{{Even still}, there is much debate over the universality of such a perspective. While this form of fitness is intuitive in the sense of evolutionary terms (``\textit{survival of the fit}"), there are cases where fitness must arguably be measured in some other way \cite{Peacock2011}, particularly when the model proposed does not explicitly consider population dynamics (model \eqref{foragingsuccess}, for example). Alternative measures exist, focusing instead on success or efficiency through some intermediate measure that many ecologists agree will correlate with fitness. These are the measures we discuss in Section \ref{measures:resources}.}

{Finally, there are cases where population dynamics nor explicit resources are included (e.g., models \eqref{scalar-eqn-1-1}-\eqref{system-eqn-3}). How, then, can we evaluate the utility of different cognitive mechanisms? This is not an easy question to answer, but \cite{Lewis2001} provides an interesting framework to measure fitness in such cases. We explore this framework in more detail in Section \ref{subsec:nopopnores}.}

Ultimately, a measure of success is highly dependent on context and should be treated as such. {For the remainder of Section \ref{sec:measuresofsuccess}, we hope to clearly describe some existing measures of success in each of the contexts outlined above. It is through {differences in} these measures that we can begin to determine the utility of different cognitive mechanisms.}

\subsubsection{Measures with Population Dynamics}\label{sec:popmeasure}

In {some formulations}, as in the models appearing at the end of Section \ref{sec:implicitmemory} and all found in Section \ref{sec:explicitmemory}, population dynamics are included. {{Then}, we may propose a number of measures of success,} where success is measured with respect to growth rates and/or survival of the population. {In this sense, the following measures are most closely related to an evaluation of the (relative) fitness of a population in terms of its ability to propagate or persist.}

The first measure, motivated by net growth rates in cell growth \cite{Karowe1989}, considers the rate of change of the total population over a prescribed time interval:
\eql{
\text{Net Growth} := \int_{t^\prime}^{t_{max}} \left( \text{Change in total population density} \right) dt.
}{netgrowthefficiency}
Thus, if $u(x,t)$ measures the population density, the change in the total population is $\frac{d}{dt} \int_\O u(x,t) dx$, and so
\eq{
\text{Net growth} := \int_{t^\prime}^{t_{max}} \frac{d}{dt} \int_\O u(x,t) dx dt = \int_{t^\prime}^{t_{max}} \int_\O \frac{\p u}{\p t} (x,t) dx dt.
}
From the discussion found in Section \ref{sec:RoT}, the population should be conserved in the absence of growth/death dynamics. Thus, all movement terms vanish {after integration} and the Net Growth {is simply}
\eql{
\text{Net Growth} := \int_{t^\prime} ^{t_{max}} \int_\O f(x,t,u) dx dt,
}{netgrowth}
where $f$ describes population growth/death. A numerical exploration of this quantity seems to be the most fruitful direction, as the formulation above leaves the dependence of Net Growth on potential knowledge-based parameters implicit. {On the other hand, if the interval is chosen such that the solution $u(x,t)$ is close to some temporally constant steady state $u^*(x)$, it may then be possible to investigate changes with respect to parameters analytically. We then propose the following open problem which can provide useful insights into the actual utility of different cognitive movement mechanisms, especially in comparison to models with no cognitive mechanisms included.}

\begin{oppr}  
What are the consequences and predictions of knowledge-based movement models in the presence of reaction terms made by the measure \eqref{netgrowth} introduced above? Are there optimal ranges for key cognition-based parameters, such as perceptual radius or rate of advection towards/away from covariates, that increase fitness in this sense? How does this measure of fitness compare to a null hypothesis model with no cognitive-based mechanisms?
\end{oppr}

{We carefully note that measure \eqref{netgrowth} is not absolute. Instead, this measure is best understood as comparing fitness between differing model formulations; the higher the net growth, the more fit the population. Note also that such a measure is influenced only by the explicit population dynamics of the model, and so the influence of different cognitive mechanisms may be lost when a zero-flux boundary condition is used.}

{Alternatively, we may use} a comparison of competitive outcomes, assuming that both survival and extinction are possibilities within the framework considered. In this sense, success is obtained when the population persists as $t \to \i$. A foraging strategy is unsuccessful if the population goes extinct. This is considered in \cite{Cho2013}, where it was shown that starvation-driven diffusion may reverse outcomes predicted by constant diffusion alone. Hence, starvation-driven diffusion may be a more successful foraging strategy. Readers should carefully note that this does not mean that starvation-driven diffusion is the more correct mechanism, but rather that it is an {alternative} explanation should we find results in nature inconsistent with predictions made by simpler models. {Such outcomes} {can} be studied analytically through changes in a principal eigenvalue or basic reproduction number, which {often} determines {the} local or global stability {of a steady state} \cite{Hess1991}. {This is a rather binary view of success, however. There are also ways to modify this perspective to determine success (or, relative success) in relation to the size or changes in the size of this key value.}
{To be concrete, consider the following
\eq{
\frac{\p u}{\p t} = D \D u + r u (1-u)
}
in a one dimensional domain $(0,L)$ subject to homogeneous Neumann boundary data. It is well known that the stability of the constant steady state $u^*=1$ is given by the sign of the principal eigenvalue, which in this case is $\mu_1 = r>0$, and so $u\to 1$ as $t \to \i$ for any $r>0$ from the theory of monotone flows \cite{Hess1991,Zhao2011dynamical}. From a biological perspective, this principal eigenvalue is exactly the net growth rate of the population, and hence can be viewed as directly connected to the fitness of the population. Taking the derivative of $\mu_1$ with respect to the parameter $r$ yields $\mu^\prime (r) = 1 > 0$, and so increasing $r$ is a quantity which contributes to the success of the species in terms of survival. In contrast, the stability does not depend on the domain size $L$ or the diffusion rate $D$, which can be observed directly since $\mu_1$ does not depend on either of these variables. On the other hand, if one considers the same problem subject to homogeneous Dirichlet boundary conditions, it is well known that the principal eigenvalue $\mu_1 = \frac{r}{D} - \left( \tfrac{\pi}{L} \right)^2$ depends on $D$, $r$, and the domain size $L$ \cite[Ch. 1.6.2]{cantrell2004spatial}. Taking the derivative of $\mu_1$ with respect to these variables yields
\eq{
\frac{d \mu_1}{d D} = -\frac{r}{D^2} < 0 , \quad \frac{d \mu_1}{d r} = \frac{1}{D} > 0 , \quad \frac{d \mu_1}{d L} = 2 \frac{\pi^2}{L^3} > 0,
}
and so we see that increasing $r$ and $L$ are beneficial strategies, while increasing $D$ is not.
}

An analysis of this sort could be applied to a wide range of models, with the difficulty dependent on the difficulty of the stability properties of the problem at hand, providing phenomenological insights into the impacts certain factors have on population {dynamics}. In the example above, the relation was found to be {monotone (either decreasing or increasing)}. More interesting examples may find non-monotone behaviour, providing an optimal \textit{window} in which certain strategies are most successful.

\begin{oppr}  
How might we formulate useful eigenvalue problems related to memory and nonlocal perception {when population dynamics are also included}? This will necessarily include a careful treatment of the nonlocal effect of the perceptual kernel, with attention directed toward the top-hat detection function. Is it possible to determine the dependence on quantities such as perceptual radius, advection speed, rates of diffusion, {or} memory uptake and decay rates?
\end{oppr}

{In contrast to measure \eqref{netgrowth}, an analysis of eigenvalues will implicitly incorporate different cognitive mechanisms considered, and so may be a better (if not more mathematically challenging) way to assess a populations' fitness. In particular, a linearization about the trivial steady state provides explicit insight into a population's fitness in terms of reproductive success.}

\subsubsection{Measures with Explicit Resources}\label{measures:resources}

{Different from models with explicit population growth/decay, many models discussed here do \textit{not} feature such processes {(e.g., model \eqref{foragingsuccess})}. {When} the total population is conserved, it does not make sense to use measures of the sort described in the previous section. Indeed, the Net Growth measure \eqref{netgrowth} would be zero!} 

{Instead, an alternative} quantitative measure of success {comes from the perspective of} {optimal }foraging {theory}. This quantity is referred to as the \textit{foraging success} and has been used as a measure to determine optimal diffusion and advection rates for model \eqref{foragingsuccess}. {In \cite{Fagan2017}}, foraging success is defined as
\eql{
\text{Foraging Success} := \int_{t^\prime} ^{t_{max}} \int_\O u(x,t) m(x,t) dx dt,
}{foragingsuccessdefn}
where the interval $(t^\prime, t_{max})$ is chosen so that the transient dynamics of the solution $u(x,t)$ have settled down. {In this context, \textit{transient} refers to the solution behaviour for smaller time intervals; this is in contrast to asymptotic dynamics, where the solution is near an equilibrium state of some form.} The purpose of this measure is to quantify the effectiveness of consumer-resource tracking, in this case, based on nonlocal information. An important distinction to make here is that this measure considers only the ability to \textit{track} the resources only, as opposed to the tracking \textit{and the utilization} of resources. As noted in \cite{Fagan2017}, however, model \eqref{foragingsuccess} does not consider mutual interference or resource depletion, and so this measure may only be valid under the assumption of sparsely populated regions and that the resources degrade more quickly than they are depleted by foragers. This may then introduce further issues, as one of the common assumptions for a knowledge-based animal movement model framework is that there is a high enough population density so that the PDE model derived is an appropriate mean-field approximation of the population density \cite{Potts2019}. 

{More recently, the foraging success as given in \eqref{foragingsuccessdefn} has been modified as a \textit{foraging efficiency} (FE) measure using the so-called \textit{Bhattacharyya coefficient} \cite{Bhattacharyya1943}, which quantifies the overlap between two distributions. In \cite{Gurarie2021}, it is given by
\eql{
\text{Foraging Efficiency} := \frac{1}{t_{max}} \int_{0}^{t_{max}} \int_\O \sqrt{u(x,t) m(x,t)} dx dt.
}{foragingefficiencydefn}
When scaled such that $\as{\O}^{-1} \int_\O u dx = 1$ and ${t_{max}} ^{-1} \int_0 ^{t_{max}} \int_\O m(x,t) dx dt = 1$, the foraging efficiency coefficient lies between $0$ and $1$ (to see this, apply H\oo lder's inequality). The restriction on the resource so that this scaling holds assumes that the average total amount of resource through, say, a year is $1$. In this sense, the Foraging Efficiency measure is normalized between $0$ and $1$, whereas the Foraging Success measure is not.}

{In the forms described above, it may be of use to note that these measures are best interpreted as a \textit{predictor} of successful foraging rather than as a \textit{result} of successful foraging, and the level of success is directly proportional to the correlation between the population density and the resource density. In other words, if a population distribution overlaps closely with the resource distribution (i.e., \eqref{foragingsuccessdefn} or \eqref{foragingefficiencydefn} is large), then the population is expected to forage successfully by virtue of being in the right place at the right time. In contrast, some authors consider the process of foraging to be ``all the methods by which an organism acquires and utilizes sources of energy and nutrients" \cite{Koy2007}. This includes key components of location, consumption, retrieval, and storage. Furthermore, some authors consider a foraging efficiency to be ``\textit{metabolizable energy gained while foraging divided by total energy spent while foraging}" \cite{Whitford2020}. This motivates a {push} for alternative measures of foraging success as a \textit{result} of successful foraging to compliment the existing correlational measures. This is not a trivial task, as such a measure would require additional mechanisms in the model to determine how much energy is used while foraging. In particular, models proposed here do not (explicitly) incorporate such mechanisms{, and hence} we use these different perspectives to propose a few open problems.
}

\begin{oppr}  
How does the analysis and predictions made in \cite{Fagan2017} for model \eqref{foragingsuccess} using the foraging success measure \eqref{foragingsuccessdefn} compare with predictions made by the foraging efficiency measure \eqref{foragingefficiencydefn}? Are they consistent in their assessment of optimal foraging strategies? For knowledge-based movement models more generally, what are the predictions made by each measure?
\end{oppr}

\begin{oppr}  
How can we formulate different measures of foraging success that incorporate mechanisms beyond mere correlation with a resource density? More precisely, how can we meaningfully quantify costs (metabolic or otherwise) associated with foraging as a means to measure a species' fitness as a result of successful foraging? 
\end{oppr}

\subsubsection{Measures with No Population Dynamics or 
 Explicit Resources}\label{subsec:nopopnores}

In many of the models introduced in Section \ref{sec:implicitmemory}, there are no population dynamics \textit{or} explicit resources! {Consequently, {previous} measures may be inappropriate, at least without further modification or justification.}

{Despite this apparent challenge, there are existing works which provide some insight into the construction of a biologically relevant fitness function in relation to adopting a particular movement strategy. In \cite{Lewis2001}, the authors study a movement model with dynamic scent marking similar to the form of model \eqref{system-eqn-2} with two wolf packs, $u$ and $v$. While {only movement is considered}, they define a geometric growth rate $R_u$ for the wolf pack $u$ as
\eq{
R_u = S \cdot N_u,
}
where $S$ is the probability that the alpha female survives the year to breed in the spring, and $N_u$ is the number of offspring that survive assuming the alpha female breeds. The quantity $R_u$ gives an expected number of offspring produced in a single year.}

{Then, the authors construct a prey density $h(x,t)$ (see \cite[Appendix A]{Lewis2001} for details) depending on the density of each wolf pack. It is worth noting that while the equation for $h(x,t)$ depends on the dynamics of each wolf pack, the movement equations do not, and so $u$, $v$ are decoupled from $h(x,t)$. This prey density $h(x,t)$ essentially acts as a resource density as in Section \ref{measures:resources}, which yields an explicit form for the expected number of yearly offspring produced:
\eq{
N_u = \s \psi \int_\O u(x) H( u(x), v(x) ) dx.
}
Here, $\s$ is a conversion rate of prey into offspring, $\psi$ is the encounter rate, and $H(u,v)$ is the average prey density during the year given space use patterns $u(x)$ and $v(x)$.}

{Finally, the authors assume that there is some probability of a wolf being killed through inter-pack aggression, and that this is proportional to an encounter rate between individuals $\a u(x) v(x)$. {Given a natural mortality rate $\mu_0 > 0$,} the overall death rate $\mu$ is given by
\eq{
\mu = \mu_0 + \a \int_\O u(x) v(x) dx,
}
The survivorship $S$ is then $e^{-\mu}$, and the fitness function for pack $u$ is given by
\eql{
r_u = - \mu_0 + \ln (\s \psi) - \a \int_\O u v dx + \ln \left( \int_\O u H(u,v) dx \right).
}{packfitnessfunction}}
{
The authors then explore optimal strategies depending on different movement mechanisms, which in their case is in differing rates of advection towards or away from covariates.
}

{
This procedure provides a route to explore a population level fitness, and in particular to compare the relative fitness between two (or more) populations within a single model. Of course, there is no reason in principle why this could not be further adapted to assess the relative fitness between two populations described by different models, assuming that the underlying features of each perspective remain roughly comparable (e.g., two separate models featuring different perceptual kernels or perceptual radii). The challenge, it seems, is instead in constructing a biologically reasonable fitness function in the style of \eqref{packfitnessfunction}, which requires some additional argumentation to justify what, exactly, the ``prey" or ``resource" is. For the case of wolf packs, some knowledge of the deer population is required; for a different species, an entirely different construction of $N_u$, or even the death rate $\mu$, would be required.
}

\begin{oppr}
Can one perform a similar analysis of optimal movement strategies with a fitness function of form \eqref{packfitnessfunction} when other cognitive movement mechanisms are included? In particular, how does this measure of pack fitness change with respect to different perceptual kernels or perceptual radii? This may be done for model \eqref{system-eqn-2} without modification; to perform a comparable analysis for other movement models without growth dynamics or explicit resources, one must first modify \eqref{packfitnessfunction} appropriately for the setting or species under consideration.
\end{oppr}

\subsubsection{Measures with Both Population Dynamics and Explicit Resources}\label{subsec:bothpopandres}

{Finally, we briefly discuss scenarios in which both population dynamics \textit{and} explicit resources are included. {O}ne could potentially choose any of the measures proposed in Section \ref{sec:popmeasure} or Section \ref{measures:resources}{;} one could also follow the framework outlined in Section \ref{subsec:nopopnores} without the need to define a resource distribution since it is already available. In fact, it could be quite interesting to explore comparisons of each of the measures proposed here for a model that includes all components necessary to make such comparisons. We highlight this as the following open problem.}

\begin{oppr}
For a movement model which includes cognitive movement mechanisms, explicit resources, and birth/death processes, how do various measures of success proposed here compare? Are they consistent with each other, or do they have a fundamental disagreement in which strategies are optimal in each sense of ``success"?
\end{oppr}

\subsection{Comparisons through Measures of Success}

We begin with model \eqref{foragingsuccess}. In \cite{Fagan2017}, the objective is to use the \textit{foraging success} \eqref{foragingsuccessdefn} as a measure to compare different diffusion and advection rates, as well as differing forms of perception kernels and perceptual radii. Interestingly, the foraging success is found to be non-monotone with respect to perceptual radius $R$ in certain environments, see \cite[Figure 3]{Fagan2017}. On the other hand, the foraging success is monotonically increasing with respect to the advection rate $\a$ under the special case when $R \to 0^+$, see \cite[Figure 2]{Fagan2017}. This provides insights into: \textit{What is the optimal detection scale for nonlocal information gathering?}; \textit{For what kinds of movement and detection does nonlocal information provide a benefit?}; or \textit{For what kind of landscapes is nonlocal information most useful?} One important aspect to note here is that these questions are possible to answer due to an available measure to infer comparatively \textit{better} or \textit{worse} scenarios within the context of this single model. 

Model \eqref{foragingsuccess} is most closely related to the new models \eqref{scalar-eqn-averagedensity} and \eqref{scalar-eqn-percapitadensity} in that the resource density is given {a priori}. There are no existing results concerning these models, and so this opens an interesting avenue for investigation.
\begin{oppr}
Which knowledge-based strategies are optimal in terms of the foraging success \eqref{foragingsuccess} for models \eqref{scalar-eqn-averagedensity} and \eqref{scalar-eqn-percapitadensity}? How do these insights compare to results obtained in \cite{Fagan2017}? It may be the case that one strategy is not always optimal, and so it becomes interesting to consider under which conditions or under which model formulation a given strategy is advantageous. This can also provide insight into how much optimal strategy predictions depend on the underlying model itself.
\end{oppr}


{For mathematical analysis, investigation of changes in a principal eigenvalue as described in Section \eqref{sec:popmeasure} may be the most fruitful direction. Ultimately, this will require the development of techniques to deal with nonlocal components which may be highly non-trivial in itself. For example, is it possible to formulate a nonlocal eigenvalue problem in a variational setting? Such developments would have high potential for interesting biological insights beyond a numerical analysis alone.}



\subsection{Emergence of Patterns}\label{sec:patterns}

One of the {core topics of investigation} in the field of movement ecology is the {connection between animal movement behaviour and the} generation of patterns in space, time or both \cite{Nathan2008}, {particularly in the absence of environmental heterogeneity \cite{Sayama2006,Potts2022b}. Indeed, countless models with explicit environmental heterogeneity have been developed, but in these cases it is much less of a surprise for non-constant steady states to emerge. Instead, it is interesting to study the emergence of non-constant steady states in problems where there is no environmental heterogeneity. This includes the formation of home ranges \cite{Moorcroft1,Moorcroft2,van2009memory} and territories \cite{Potts2016, Potts2016-2, Potts2019} and the potential mechanisms required for each.} Classically, this is discussed through linearization techniques which gives rise to so-called \textit{Turing instability} \cite[Ch. 2]{Murray2003}, where a constant steady state can be destabilized {by} the presence of diffusion. There are now a significant number of resources exploring in great detail the possible outcomes {for classical reaction-diffusion equations}, including not only \textit{when} patterns will form, but also which \textit{types} of patterns one might expect to find. While the original motivation {of advection equations such as \eqref{prototype}} was to {describe} cell growth phenomena, some of these same ideas and motivations transfer to the study of movement ecology.  

As it currently stands, we have little understanding of when patterns may emerge, let alone which types of patterns we might expect to appear {for non-standard PDE models explored here}. For the implicit memory models discussed in Section \ref{sec:implicitmemory}, the same linearization techniques used in classical setting yield some preliminary insights: foragers must have sufficient uptake of memory, and cannot forget experiences too quickly in order for patterns to emerge; furthermore, the advection rate must be sufficiently large compared to the diffusion rate in order for patterns to form \cite[Section 3.1]{Potts2016} (see also \cite[Section 4]{Potts2019}). However, we still do not have a robust understanding of necessary and sufficient conditions for patterns to emerge, {hence the appeal to a numerical exploration of local stability properties}. A similar trend holds for the explicit memory models discussed in Section \ref{sec:explicitmemory}: the advection rate must be sufficiently large in relation to the diffusion rate for constant steady states to {become} destabilized (see, e.g., \cite[Corollary 3.9]{Shi_JDDE}, \cite[Theorem 2.3]{Song2021}, \cite[Theorem 3.3]{Shi2021a}). Of course, the techniques used in these two settings are quite different.

From a mathematical perspective, this yields a rich area of analysis which can be roughly understood across all models {introduced thus far}: when do patterns form, and what kinds of patterns are possible? It is not clear what kinds of patterns are possible for implicit memory models; {most} existing works consider only the possibility of patterns emerging in relation to parameters appearing in the model, e.g., \cite{Potts2016, Potts2016-2, Potts2019}; {in some cases, effort has been made to describe possible non-constant steady states using energy methods in {the local case}, see \cite[Section 3.4]{Potts2016}.} We conjecture that some models may in fact have infinitely many piece-wise constant steady states {in this local limit}. {Recent progress has been made in \cite{Giunta2022}, where the authors rule out certain pattern forms (see Proposition 2 and Remark 2) and verify the existence of patterns that are {either} roughly piecewise constant or spikey solutions.} In the case of explicit memory models, we have a better understanding of what types of patterns may form, at least from a simulation perspective: spatially constant temporally periodic steady states, temporally constant spatially varying steady states (stripe patterns), and spatiotemporally varying steady states (checkerboard patterns) have all been observed, see \cite{Song_JDE, Shi_JDDE, Song2021, Shi2021a}. However, it is still not clear exactly when these patterns are expected to form, and whether these are all of the possibilities (we believe not).

Ultimately, linearization techniques alone are not enough to obtain the deepest insights. While they provide insights into when patterns may form, these techniques do not tell us which patterns will emerge, and more importantly do not tell us the mechanisms behind specific patterns that emerge. {This leads to another broad open problem applicable to virtually any model introduced here.}

\begin{oppr}  
How do new modelling components introduced in this review affect the stability of constant steady states? Can we identify and classify the possible patterns? What are the mechanism behind the generation of certain patterns? Can we connect existing simulations with analytical insights? Can these be connected to patterns found in nature, where data is available? \cite{Potts2020} provides some insight in regards to connections with data. Finally, how should we interpret these patterns, both biologically and mathematically? All patterns are ultimately mysterious: are they typical? Are these strange cases or normal for these models?
\end{oppr}

\subsection{Travelling Wave Solutions}\label{sec:travwave}

{{For problems} in an unbounded domain, a popular area of application in biology and other sciences is the idea of \textit{travelling wave solutions}. In the setting of mathematical ecology, these solutions are often interpreted as the front of a population density curve, indicating expansion or retraction of the population. To fully explore the wide range of works studying such solutions is far beyond the scope of this review. Instead, we make reference to the rather comprehensive review paper by Volpert and Petrovskii \cite{Volpert2009} which covers {many} biological applications of travelling waves. Readers are also directed to \cite[Ch. 13]{Murray2003} or \cite[Ch. 4]{Perthame2015} for a more elementary introduction and additional references.}

{The main idea is to seek solutions of the form $u(x,t) = w(x - ct)$, where $c$ is a constant wave speed which, in general, depends on the form of the equation considered. Many problems are formulated for a scalar reaction-diffusion equation in domain $\O = \mathbb{R}$. In the current setting, {we are interested in} the properties of travelling wave solutions {as they relate to} cognitive mechanisms. Based on the discussion in Section \ref{sec:considerations} concerning the definition of perceptual kernels near the boundary and how to appropriately choose boundary conditions, the travelling wave setting offers a potentially fruitful area of investigation since no further modification of the kernel is required.}

{Once a model has been formulated in a framework appropriate for analysis, there are a handful of key questions of particular interest in biological and ecological application: existence of travelling waves, stability of travelling waves, and the speed of propagation. Other questions include the effect of the initial condition and rates of convergence to a travelling wave solution. Existence of a travelling wave solution, in a rigorous mathematical sense, is of interest for the same reasoning described in Section \ref{sec:considerations}: there are cases where numerical simulation will produce {something resembling} a travelling wave solution while it is known that no such solution exists \cite[Ch. 1]{petrovskii2005exactly}! Stability properties are also of interest, which is generally determined in a similar fashion to reaction-diffusion equations in a bounded domain: linearization about the constant steady state and its spectrum. However, travelling wave solutions generate a family of solutions once one is obtained, and so some care must be taken. Calculation or estimation of the wave speed is also of particular interest since this describes the rate at which an invasion wave will move. Depending on the context, it may be of interest to slow down or speed up this wavefront, and so an explicit calculation of the wave speed as it depends on model parameters is ideal. In many cases, though, an explicit formula is not possible and estimation may be the only option.}

{As far as we are aware, there are no existing works which consider such questions in models that feature cognitive mechanisms. This is almost certainly due to the recency of models described here, but also because of the technical difficulties due to nonlocal terms appearing at a higher order. As such, we formulate a handful of open problems in this direction.
}

\begin{oppr}\label{nonlocaltravel}
Do travelling wave solutions exist solving the nonlocal problem \eqref{scalar-eqn-1-1} when the detection function is chosen to be any of the forms defined in \ref{detectionkernels}? In cases where the answer is affirmative, is it possible to then compute or estimate the wave speed depending on the form of detection function or perceptual radius? What role do the rates of advection and/or diffusion play?
\end{oppr}

\begin{oppr}
Do travelling wave solutions exist solving the discrete delay problem \eqref{delayprototype} for fixed delay $\tau>0$? This can be explored for either $f(u) \equiv 0$ or with a general birth/death process $f(u)$. In cases where the answer is affirmative, is it possible to calculate or estimate the wave speed? How does the wave speed depend on the delay parameter $\tau$ or the rates of advection and/or diffusion? We remark that while travelling wave solutions have been shown to exist for models with delay at the lowest order, the challenge here is the fact that the delay parameter appears underneath the gradient.
\end{oppr}

\begin{oppr}\label{disttravel}
Do travelling wave solutions exist solving the distributed delay problem \eqref{delayfinal} for temporal kernels of the form defined in \eqref{weakstrongkernel}? In cases where the answer is affirmative, is it possible to compute or estimate the wave speed? How does the wave speed depend on the delay parameter $\tau$ in the temporal kernel or the rates of advection/diffusion?
\end{oppr}

{In each of the open problems described above, only scalar models are chosen as they are almost certainly easier to solve compared to systems. This may introduce further challenges for Open Problem \ref{disttravel}, since the most promising technique seems to be converting the distributed delay problem into a system of two or more equations depending on the temporal kernel chosen. On the other hand, there also exist many works exploring the existence of travelling wave solutions to systems. This motivates the following open questions.
}

\begin{oppr}
Do travelling wave solutions exist solving the nonlocal system \eqref{system-eqn-1} when the detection function is chosen to be any of the forms defined in \ref{detectionkernels}? In cases where the answer is affirmative, is it possible to then compute or estimate the wave speed depending on the form of detection function or perceptual radius? What role do the rates of advection and/or diffusion play? Resolving such a question would be a direct extension of Open Problem \ref{nonlocaltravel}.
\end{oppr}

\begin{oppr}\label{nonlocalsystemtravel}
Do travelling wave solutions exist solving the nonlocal system \eqref{system-eqn-2} for a single species (i.e., a single species with a dynamic cognitive map) when the detection function is chosen to be any of the forms defined in \ref{detectionkernels}? In cases where the answer is affirmative, is it possible to compute or estimate the wave speed? How does the wave speed depend on detection function, perceptual radius, or other parameters appearing in the equation? Can these results be extended to $n$ interacting species with cognitive map? 
\end{oppr}

{We conclude by again emphasizing that the mathematical challenge here is the presence of nonlocality at a higher order. Of interest may be the work of Wang \cite{Wang2013}, which explores the mathematical features of travelling waves for a class of chemotaxis models. While these forms do not incorporate nonlocal terms, it may at least provide insight into possible techniques to apply, particularly for Open Problem \ref{nonlocalsystemtravel}.}

\subsection{Critical Domain Sizes}\label{sec:critdom}

{An important area of theoretical spatial ecology is the concept of \textit{critical domain size} or \textit{critical patch size}. In one spatial dimension, this is a phenomenon where, under certain conditions, a problem may have a critical domain length $L^*$ such that the population goes extinct {(persists)} when $L < L^*$ {($L > L^*$)} {This} is of particular interest as it gives insight into the relationship between a population's persistence and the size of its habitat. This can be generalized to higher dimensions, though it can become much more complicated due to regularity or other special geometric properties of the domain. Readers are directed to \cite{cantrell2004spatial} or \cite{Murray2003} for an introduction to these problems in {several} contexts.}

{W}e must first consider {which} mechanisms introduce such a threshold. For models {with a conserved} population, it is relatively trivial to note that a critical domain size cannot exist. Instead, such a critical threshold depends heavily on the boundary conditions chosen, and in particular on conditions that result in a loss of population across the boundary. This can be observed in the brief analysis done near the end of Section \ref{sec:popmeasure}: in the case of a Neumann boundary condition, no population is lost across the boundary, and the positive steady state is always stable; this is in contrast to the homogeneous Dirichlet boundary condition, where the stability of the positive steady state now depends on the domain size $L$.

In the context of cognitive mechanisms, a potentially substantial difficulty {is introduced}: we must consider boundary conditions that require modification of the detection function near the boundary. Since this has not been done before, there is the auxiliary step of understanding a given equation with nonlocal components in a bounded domain. Hence, we introduce the following open problem.

\begin{oppr}
Consider problem \eqref{scalar-eqn-1-1} in a bounded, one dimensional spatial domain $(0,L)$ subject to homogeneous Dirichlet boundary conditions and a top-hat detection function defined in the style of \eqref{cutoffkernel}. Does there exist a critical threshold $L^*$ such that the only steady-state solution is the trivial one for $L < L^*$, and a positive solution exists for $L>L^*$? If this is proven to be the case, how does the critical value $L^*$ depend on the perceptual radius, or advection and/or diffusion rates? Can this be generalized to the $n$-species model \eqref{system-eqn-1}?
\end{oppr}

\subsection{Numerical Analysis}\label{sec:numericalanalysis}

{We conclude this section with a brief discussion of numerical techniques  as they apply to models introduced in Section \ref{sec:implicitdynamicmemory}. In general, it is challenging to deal with nonlocal terms in a numerical setting in a bounded domain {due to the nonlocality of the equation}. As such, in-built PDE solvers (such as MATLAB's ``pdepe" function) cannot solve these problems even in one spatial dimension, and so new solvers must be built. {T}here is at least one exception to this rule, and it happens to be the same exception as in the analytical discussion: a periodic boundary condition does not require further modification of the detection function. For this reason, {so-called} \textit{pseudo-spectral} methods offer a very promising route to efficiently and quickly simulate such solutions.
}

{The main idea behind pseudo-spectral methods is to solve the problem in space using Fourier series, while solving the problem in time using some other standard scheme (e.g., forward Euler or Runge-Kutta methods). The biggest advantage {here} is the convolution theorem: the Fourier transform of a convolution of two functions is equal to the product of the Fourier transform of each function. In mathematical terms, $\mathcal{F} (f \star g) = \mathcal{F} (f) \mathcal{F}(g)$, where $\mathcal{F}$ denotes the Fourier transform {and $\star$ denotes convolution}. In many cases it is possible to compute the Fourier coefficients for a given detection function. As such, a pseudo-spectral method essentially removes the difficulty that the spatial convolution {can bring in finite-difference schemes}. We direct readers to \cite{Giunta2021} which discusses this technique in more detail. We leave this section with the following general open problem.}

\begin{oppr}
Is there a general technique to efficiently solve nonlocal problems of the form appearing in Section \ref{sec:implicitdynamicmemory} subject to boundary conditions that are not of the periodic type? 
\end{oppr}

\section{Mathematical Techniques and Challenges}\label{sec:mathtech}

\subsection{Well-posedness of Existing Models}\label{sec:wellposedness}

\subsubsection{Implicit Memory Models}\label{sec:wellposedimplicit}

We begin with a discussion of some fundamental results (or lack thereof) in the implicit memory models {found} in Section \ref{sec:implicitmemory}. {Many open problems presented in the aforementioned section provide rather explicit questions currently without resolution.} Given the relatively sparse literature available {for} these systems, it is easiest to describe the cases where existence is \textit{known}, {complementing {previous} open problems}. 

In general, model \eqref{prototype} has a large body of literature concerning the well-posedness of the problem under appropriate regularity assumptions on the function $a(x,t)$. If $a(x,t)$ is twice continuously differentiable, for example, a unique solution exists for all time. Thus, any model which provides $a(x,t)$ explicitly enjoys well-posedness, essentially following from standard theory of linear partial differential equations, see e.g., \cite[Ch. 8]{Wu2006Elliptic}, \cite{friedman1964partial}, \cite[Part II]{evans1998partial}. This means that models \eqref{foragingsuccess}, \eqref{scalar-eqn-densite}-\eqref{scalar-eqn-averagedensity} have a unique, global solution, so long as the resource $m(x,t)$ is sufficiently smooth. On the other hand, model \eqref{scalar-eqn-percapitadensity} is nonlinear and existence is not so obvious. 

The difficulty is increased substantially with models introduced in Section \ref{sec:implicitdynamicmemory}. Indeed, such systems are nonlinear at a higher order (i.e., the coupling occurs inside the gradient). Proving the well-posedness of such models is highly nontrivial, with at least one exception: if the perceptual kernel $g(\cdot)$ is sufficiently smooth, the regularity requirements for the nonlinear term can be transferred to the kernel itself, and so existence follows from standard $L^p$-theory of parabolic equations \cite[Ch. 9]{Wu2006Elliptic}. Recently, the well-posedness of model \eqref{system-eqn-1} has been shown under the assumption that the perceptual kernel is twice continuously differentiable \cite{Giunta2021}. This was shown using semi-group theory {and fixed point methods} \cite{Amann1995, Pazy1983}. We note that a similar existence result should be possible using $L^p$-theory of parabolic equations when the perceptual kernel is twice differentiable, and its second derivative belongs to $L^\i (Q_T)$,{ where $Q_T = \O \times (0,T)$}. Indeed, estimates found in \cite[Ch. 9]{Wu2006Elliptic} hold somewhat trivially, since we may transfer the derivatives from the solution $u$ itself to the spatial kernel instead. Hence, the existence of a strong solution belonging to $W^{2,1}_2 (Q_T)$ for a smooth, bounded domain $\O$ and $T>0$ fixed follows from a mere $L^2(Q_T)$ estimate on the solution $u$. {Of course, this is assuming one has determined an appropriate method to deal with whatever boundary conditions are being applied. As previously noted, periodic boundary conditions are easiest to deal with in this sense. In fact, this is roughly the technique used in \cite{Jungel2022}, where the authors study the existence and uniqueness of solutions to a special case of model \eqref{system-eqn-1}.} 

As a summary, the well-posedness of any model featuring a Gaussian kernel enjoys well-posedness due to the regularity of the kernel, {but one is restricted to studying the problem in an infinite domain}. The exponential kernel is more difficult due to the lack of differentiability at zero. The most difficult case, and perhaps most biologically interesting, is the top-hat detection function, which is discontinuous. This fact is precisely what makes these models so difficult to study analytically, and new techniques or tools must be used. On the other hand, modellers may be satisfied with approximating the top-hat detection function through a regularization technique, e.g., mollification, in which case well-posedness is less of an issue and will follow from a careful application of existing theory.

\subsubsection{Explicit Memory Models}\label{sec:wellposedexplicit}

In contrast, {time-delay models} found in Section \ref{sec:explicitmemory} generally have more results concerning their well-posedness. {In fact, almost no modification of standard techniques is required to prove existence of a solution.} In this sense, their complexity is offset by the existing literature for delay equations and do not often require a study of nonlocal (in space) effects. We highlight some key results existing in the literature:
\begin{itemize}
\item Model \eqref{delayprototype} has a unique, global solution under mild conditions, see \cite[Proposition 2.1]{Shi_JDDE}. This follows from a standard bootstrapping method, assuming regular initial data and applying classical theory of parabolic equations.
\item \cite{Song_JDE} does not consider the existence of solutions to model \eqref{delayprototypenon-local} (a proposed open problem). This may follow from standard arguments, since the non-locality appears in the reaction term rather than inside the gradient.
\item Model \eqref{delayprototypenon-localtime} has a unique, global solution under mild conditions, see \cite[Theorem 2.1]{Shi_Nonlinearity}. This also follows from a standard bootstrapping argument with minor modifications to deal with two time delays.
\item \cite{An_DCDS} considers the existence of inhomogeneous steady states to model \ref{delayprototypenon-localtime2} under some stronger regularity assumptions on the growth term and nonlocal kernel appearing within it. This is achieved through a Lyapunov-Schmidt reduction, see \cite[Theorem 2.1]{An_DCDS}. Existence of solutions to the time dependent problem is not considered, and so this is an open {problem}.
\item \cite{Song2021Preprint} does not consider the existence of solutions to problem \eqref{delaymodelsystem}, and so this is an open {problem}.
\item\cite{Shi_JMB} does not explicitly consider the existence of solutions to model \eqref{delayfinal}, however \cite[Lemma 2.2]{Shi_JMB} (the statement of which is found in the Appendix) provides an interesting equivalence result between the delay differential equations and a certain form of Keller-Segel chemotaxis model. Hence, existence and non-existence of solutions may follow from the vast literature concerning Keller-Segel type systems.

\item \cite{Song2021} does not consider the existence of solutions to model \eqref{delayprototypenon-localtime-dist}, however some references concerning periodic solutions and travelling wave solutions are provided. The study of the well-posedness of this model is an open {problem}.
\end{itemize}

To conclude, models which feature a spatiotemporal convolution, defined most generally in line \eqref{percepspatim}, do not have many results concerning the well-posedness of solutions. Similar to the caveats made in Section \ref{sec:wellposedimplicit}, existence is less of an issue when the kernels are appropriately smooth. The most interesting cases, however, do not often feature such regularity (as in the case of the top-hat detection function {due to a lack of continuity}, or when the kernel is taken to be the fundamental solution to the heat equation {and the kernel becomes singular}), and thus further development of analytical tools is necessary.

\subsection{Pattern Formation through Linear Stability Analysis}\label{sec:linstab}

\subsubsection{Implicit Memory Models} 
We discuss some of the results found in \cite{Potts2016} and the general trends based on dispersal relations (see \cite[Figure 1]{Potts2016}), at least as they are found in the case of two interacting species. These results are most closely related to model \eqref{system-eqn-3}, however readers are reminded that the model appearing in \cite{Potts2016} treats the memory component $k_i$ as a probability distribution as opposed to a magnitude alone (see the brief discussion in Section \ref{sec:implicitmemory} before model \eqref{system-eqn-3} is introduced). 

First, it is found that in the limit as the perceptual radius $R$ decreases, the set of wavenumbers at which patterns may form increases in size{;} in the limit as $R \to 0^+$, patterns may form at arbitrarily high wavenumbers, and so the problem becomes ill-posed. This indicates that a decreasing perceptual radius has a destabilizing effect on the constant steady states.  

Another set of important parameters considered is the rate at which memories decay ($\mu$) or the rate at which foragers update their cognitive map should they revisit a site and deem it safe ($\b$). It was found that if either $\mu$ is too large, or $\b$ is too small, patterns cannot occur at any wavenumber. Hence, for patterns to form, foragers cannot forget information too quickly, and they must have some mechanism to update their cognitive map. In the case of conflict zones, this means that foragers must have some mechanism through which they feel safe should they revisit a site and experience no conflict. Similarly, it was found that the set of wavenumbers at which patterns may form increases as the advection rate increases. This suggests that foragers must move quickly enough towards safe areas or away from conflict zones for patterns to occur. On the other hand, it was found that the rate of diffusion and rate at which conflicts occur does not change the set of wavenumbers at which patterns may form, however the rate at which these patterns grow is smaller when the conflict rate $\rho$ is larger or the diffusion rate is smaller.

This exploration provides some preliminary insights into the relative effects these parameters have on the possibility of pattern formation, and some further investigation of these effects is explored in \cite{Potts2019}, where some simplifying assumptions allow the authors to unify models \eqref{system-eqn-1}, \eqref{system-eqn-2}, and \eqref{system-eqn-3}. These simplifications are convenient for their analysis, however it may be interesting to investigate the effects that differing parameter regimes may have without {simplification}. In particular, one may consider the effects of different rates of diffusion, rates of advection, as well as the relative effects of the memory decay rates, without an appeal to a quasi steady-state approximation. Readers should note that {for} three or more populations, less is known, and the numerical results found in \cite{Potts2019} suggest that the dynamics may be chaotic.

Finally, readers should note that the above holds when we know the existence of a positive steady state. This is easiest when constant steady states are {valid}, {but} this is highly dependent on the boundary conditions chosen: {with} a homogeneous Dirichlet boundary condition, {the constant steady state is the trivial one}, and so an additional step {of proving the existence of a non-trivial steady state is required.}

\begin{oppr}  
How do different parameters change the broad insights obtained above? That is, what role might unequal diffusion rates, advection rates, perceptual radii, perceptual kernels etc. play in altering space use outcomes? Do the general trends highlighted above remain true? If not, under what conditions do these general trends change?
\end{oppr}

\begin{oppr} 
Can we expand our understanding of the impact of numerous (3 or more) interacting populations in the current setting? This is considered briefly in \cite{Potts2019}, for example, where 3 interacting populations demonstrate oscillatory behaviour. In general, increasing the number of populations will complicate the analysis significantly, and so this is a highly non-trivial question to answer in any universal sense.
\end{oppr}

\subsubsection{Explicit Memory Models} 
In the case of explicit memory models ({Section \ref{sec:explicitmemory}}), a linear stability analysis {still} allows one to investigate the possibility of pattern formation. However, due to the complex nature of delay differential equations, the resultant analysis of local stability of constant steady states is {often} significantly more involved. A good introductory reference for partial delay differential equations is found in \cite{wu1996theory}, however, readers should carefully note that all delay parameters appear in lower order terms. For this reason, new tools and techniques need to be developed in the case of knowledge based movement models since the delay parameters appear at a higher order, increasing the difficulty substantially. Despite this, we may still discuss some of the key insights found in the use of temporal delays, which is primarily done through changes in stability of possible constant steady states (model \eqref{delayprototypenon-localtime2} is the one exception due to the hostile boundary condition). Readers should note that these are general trends, and to cover every possible outcome in detail here is more challenging than in the cases without delays since the possible outcomes are much richer and can vary significantly across models.

We begin with the prototypical delay model \eqref{delayprototype}. It was shown in \cite{Shi_JDDE} that the stability of the constant steady state depends on the ratio of diffusion rate and advection speed, but is independent of the discrete the time delay $\tau$. Roughly, the advection speed away (or towards) high density areas must be sufficiently large in relation to the diffusion in order to destabilize the constant steady state. This suggests that the average time which the foragers reference back to does not influence the emergence of patterns, but there must be a mechanism by which the foragers move towards these preferred areas more quickly than the random diffusive movements. Intuitively, this makes sense: if the random motion is too large, this overtakes any possibility of aggregation/segregation and the population density will remain uniformly distributed {in} the environment. 

Model \eqref{delayprototypenon-local} generalizes this through an inclusion of nonlocal effect in the growth of the population. Interestingly, \cite{Song_JDE} shows that the changes in stability remain roughly the same, however these changes are no longer independent of the delay parameter $\tau$. Similarly, the rate of advection must be sufficiently large in magnitude {to destabilize} the constant steady state. This is shown through the appearance of a Turing-Hopf or double Hopf bifurcation for some values of $\tau$ with respect to the advection speed. This is most easily viewed through \cite[Figure 5]{Song_JDE} where a stability region is provided with respect to advection and delay parameters. Furthermore, there are many different forms of steady states that are possible (periodic solutions, non-constant steady states), as opposed to only a constant steady state. This demonstrates that a nonlocal effect in the growth term promotes a wider variety of potential outcomes in animal space use.

Model \eqref{delayprototypenon-localtime} {studied in \cite{Shi_Nonlinearity}} is most similar to {model \eqref{delayprototype} with similar conclusions drawn:} the advection rate must be sufficiently large in comparison to the diffusion rate to destabilize the constant steady state, and this occurs independent of both delay parameters (see \cite[Theorem 2.6, Theorem 2.8]{Shi_Nonlinearity} and the Discussion). However, when the advection speed is small in relation to the diffusion rate, the memory delay parameter plays a key role in determining the stability of the constant steady state.

Model \eqref{delayprototypenon-localtime2} {studied in \cite{An_DCDS}} is an exception to other delay model results, as the hostile boundary condition implies that $0$ is the only constant steady state. Therefore, {they} ensure that a non-trivial steady state exists (see \cite[Theorem 2.1]{An_DCDS}). {Then}, a similar trend holds: destabilization of this steady state is induced by a sufficiently large advection rate (see \cite[Theorem 3.9]{An_DCDS}). Numerically, it is observed that the hostile boundary condition results in ``stripe" patterns, as opposed to ``checkerboard" patterns as found in the case of zero-flux {or periodic boundary conditions.}

Different from the models discussed above, model \eqref{delaymodelsystem} {in \cite{Song2021}} features consumers as well as resources, and hence it is fundamentally different from a single species model. {T}his model features {both} memory {and} dynamic resources, which may be more realistic than resource densities {given} a priori. Despite this, some of the general trends still hold: when the advection speed is small, there are no changes in stability of the constant steady state; when the advection rate is moderate, the delay parameter may have a destabilizing effect; when the advection speed is large, there is a critical delay value $\tau^*$ for which the constant steady state is stable when $\tau < \tau^*$, and unstable when $\tau > \tau^*$ (see \cite[Theorems 2.1-2.3]{Song2021}). {Similar to maturation delays}, interaction between consumers and resources introduces a key difference from single species models in that the delay parameter can drastically change the long term dynamics.

The final two models {instead} feature distributed delays. Model \eqref{delayfinal} {studied in \cite{Shi_JMB}} features a single species model with distributed delay, as well as differing kernels determining how previous information influences movement. First, it is shown that the choices of kernels found in \cite{Song_JDE} correspond to an equivalence of systems, particularly a well know Keller-Segel system (see the Appendix). This means that the analytical tools used to study the dynamics of Keller-Segel models can be used on this distributed delay model. Despite the departure from the discrete delay models, the same trend holds once more: the advection speed must be sufficiently large in order to destabilize the constant steady state (see \cite[Theorems 3.2, 4.3]{Shi_JMB}).

Finally, model \eqref{delayprototypenon-localtime-dist} {studied in \cite{Song2021}} includes a distributed delay in both memory {and} growth term{s}. {T}he destabilizing effect of advection speeds no longer holds, at least in the weak kernel case: when the advection speed is positive (segregation effect), patterns cannot occur; when the advection speed is negative (aggregation), bifurcations can occur at any advection speed, depending on the delay parameter $\tau$ (see \cite[Theorems 2.1-2.3]{Song2021}). {T}his suggests that large advection speeds \textit{in magnitude} is no longer sufficient to induce destabilization. Similarly, when the memory delay is held fixed, bifurcations occur in relation to the advection speed and maturation delay parameter $\s$. 

\subsection{Unification and Equivalence of Existing Models}

{T}o combine these models into a cohesive form, one may generalize the spatial kernel introduced in \eqref{resourcedetection} to include both spatial and temporal influences {as in \eqref{percepspatim}}. To this end, suppose we are given a {potential} $a(x,t)$. We then define
\eql{
\overline{a} (x,t) = {\mathcal{K} * * a (x,t)} := \int_{-\i}^ t \int_\O \mathcal{K} (x,y,t,s) a(y,s) dy ds,
}{spat-time-kernel-coupled}
where $\mathcal{K}$ is some reasonably defined space-time kernel describing modifications to the quantity $a(x,t)$ with respect to both distance an time. {For example}, we may consider space-time kernels that are separable in their variables, i.e., $\mathcal{K}(x,y,t,s) = g(x,y) \mathcal{G} (t,s)${, so that}
\eql{
\overline{a} (x,t) =   \int_{-\i}^{t} \int_\O g(x-y) \mathcal{G}(t-s) a (y,s) dy ds .
}{spat-time-kernel}
{Since space and temporal effects are now independent}, one may consider two ``averaging" processes {over each domain} so that combining both yield the form above. That is, define the linear operators $S$ and $T$ by
\eq{
S [ a ] (x,t) &:= \int_\O g(x-y) a (y, t) dy , \\
T [ a ] (x,t) &:= \int_{- \i} ^t \mathcal{G} (t-s) a (x,s) ds .
}
Assuming these integrals are well defined, we may take the composition
\eq{
T [ S [ a ] ] (x,t) &= T \left[ \int_\O g(x-y) a (y, t) dy \right] = \int_{- \i} ^t \left( \int_\O g(x-y) a (y, s) dy \right) \mathcal{G}(t-s) ds = \overline{a} (x,t) ,
}
to recover \eqref{spat-time-kernel}. One may then relate all forms presented so far through various choices in $\mathcal{K}(\cdot)$, $g(\cdot)$, $\mathcal{G}(\cdot)$, as well as the quantity $a(x,t)$ itself. For example, any model that does not feature an integration of information over previous times, {we fix} $g(t) = \d (t)$. In doing so, \eqref{spat-time-kernel} collapses to the form in \eqref{resourcedetection}, and we recover the models found in Section \ref{sec:implicitmemory}. When a discrete delay is included {with no perception} (as in Section \ref{sec:discretedelays}), we {fix} $\mathcal{G}(t) = \d (t- \tau)$ and $g(x) = \d(x)$. In the case of distributed delays, one must consider the more general form \eqref{spat-time-kernel-coupled} and choose $\mathcal{K}$ to be a product of the Green's function for the heat equation and the weak or strong kernel defined in \eqref{weakstrongkernel}, which recovers models \eqref{delayfinal} and \eqref{delayprototypenon-localtime-dist}. Table \ref{table1} highlights these connections explicitly for models introduced in this manuscript. In it, the advective potential $a(x,t)$ as found in \eqref{spat-time-kernel-coupled} is described, along with the relevant kernels $g(x)$ and $\mathcal{G}(t)$, and the associated reference(s) when applicable.

\begin{table}[b!]
\caption{Models, Potentials \& Kernels}
\begin{tabular}{|p{2cm}|p{5cm}|p{5cm}|p{2.5cm}|}
\hline
Model                      & Advective Potential           & Kernel $\mathcal{K} (x,y,t,s)$          & Reference             \\ \hline
$\eqref{foragingsuccess}$  &    Resource $m(x,t)$ &    top-hat, exponential, Gaussian     &  $\cite{Fagan2017}$           \\ \hline
$\eqref{scalar-eqn-densite}$&   Euclidean distance $\norm{x-x_0}$ from den site located at $x_0$   &    $\d(x) \d(t)$ (local information only) &  $\cite{Moorcroft1,Johnston2019}$ \\ \hline
$\eqref{scalar-eqn-averagedensity}$ &  Ratio of resource density and average resource density, $m(x,t) / \overline{m}$ &   $\d(x)$, top-hat, exponential, Gaussian & Newly proposed          \\ \hline
 $\eqref{scalar-eqn-percapitadensity}$ &  Per capita resource density, $m(x,t)/u(x,t)$ &    $\d(x)$, top-hat, exponential, Gaussian  &   Newly proposed        \\ \hline
$\eqref{system-eqn-1}$  & Density of all other species $u_i$, $i=1,\ldots,n$  & top-hat & $\cite{Potts2019,Hillen2009}$\\ \hline
$\eqref{system-eqn-2}$ & Density of scent marks left on landscape by other populations, $p^i(x,t)$ &    top-hat                  &$\cite{Lewis1993, Potts2019}$\\ \hline
$\eqref{system-eqn-3}$ & Cognitive map given by memory of conflict zones, $k^i(x,t)$ & top-hat & $\cite{Potts2016,Potts2019}$ \\ \hline
$\eqref{system-eqn-RC1-case1}$ &Density of prey, $v(x,t)$; extensions include density of marks left on landscape $\eqref{system-eqn-RC1-case2}$ or memory of areas where resources were previously found \eqref{system-eqn-RC1-case3} &  top-hat &$\cite{Song2021Preprint}$ \\ \hline
$\eqref{system-eqn-3-mod1}$-$\eqref{system-eqn-3-mod2}$ & Extensions of $\eqref{system-eqn-3}$, including memory smearing and/or nonlocal conflict zones & top-hat & Newly proposed \\ \hline
 $\eqref{system-eqn-RC1}$-$\eqref{system-eqn-RC1-case3}$ & Extensions of $\eqref{system-eqn-RC1-case1}$, including different cognitive map formulations & top-hat & Newly proposed \\ \hline
$\eqref{delayprototype}$-$\eqref{delayprototypenon-localtime2}$& Population density $u(x,t)$ & $\d(x)\d(t-\tau)$ for $\tau > 0$& $\cite{Shi_JDDE, Song_JDE, Shi_Nonlinearity, An_DCDS}$\\ \hline
$\eqref{delaymodelsystem}$& Density of prey, $v(x,t)$  & $\d(x)\d(t-\tau)$ for $\tau > 0$  & $\cite{Song2021}$  \\ \hline
$\eqref{delaymodelcompsystem}$ &  Density of own population and competitors, $u(x,t)$ and $v(x,t)$ &  $\d(x)\d(t-\tau)$ for $\tau > 0$ & $\cite{Shi2021a}$    \\ \hline
$\eqref{delayfinal}$\&$\eqref{delayprototypenon-localtime-dist}$& Density of own population $u(x,t)$ &  Green's function for the heat equation, $G(x,y,t) \mathcal{G}(s,t)$, where $G$ is the Green's function for the heat equation and $\mathcal{G}$ is either the weak or strong kernel, see $\eqref{weakstrongkernel}$  & $\cite{Shi_JMB, Song2021}$ \\ \hline
$\eqref{shortlongequation}$  & Cognitive map with both short- and long-term memory mechanisms, $a(x,t)$  &  top-hat, exponential, Gaussian, general form as described in Section $\ref{sec:perception}$                     & Newly proposed \\ \hline
$\eqref{SDA-densite-1}$ & $\tilde{\omega}(s) m {+} \g \norm{x-x_0}$, starvation-driven advection and densite & $\d(x) \d(t)$ (local information only) & Newly proposed \\ \hline
\end{tabular}
\label{table1}
\end{table}

It should be noted that these generalizations do not necessarily lend themselves to a concrete mathematical analysis{;} for some cases it may provide a stepping stone to results concerning the existence and uniqueness of solutions to the time-dependent problems.

\section{Concluding Remarks}\label{sec:concremarks}

We conclude our {review} with some overarching themes and broad impacts of the works discussed {thus far}. First, we have taken care to introduce multiple key cognitive mechanisms in these pioneering models {from the biological perspective}. {T}his includes {mechanisms of} perception, memory, and learning. {Mathematically,} these are incorporated through an advective {potential} which biases movement {beyond} passive diffusion. {The details found in Section \ref{sec:derivations}, while certainly not new to this area of study, provide useful insights into the mechanistic derivation of diffusion-advection equations in general, which in turn produces a foundation on which subsequent models can stand. We then systematically explore how these cognitive modelling components commonly appear in the existing literature.}

Perception {(Section \ref{sec:perception})}, included through a spatial convolution \eqref{generaldetection}, incorporates differing perceptual capabilities through a perceptual radius $R$ and perceptual kernel $g(\cdot)$. Memory, included as an implicit static quantity {(Section \ref{sec:staticmemory})}, an implicit dynamic quantity {(Section \ref{sec:implicitdynamicmemory})}, or an explicit quantity through time delays {(Section \ref{sec:explicitmemory})}, incorporates the process of encoding, storing, and retrieving information within the equation(s) describing movement. Often, this is included through a cognitive map. Learning {(Section \ref{sec:learning})}, either implicitly through memory uptake functions, or explicitly through variable diffusion/advection rates via satisfaction measures, allows one to consider, compare and contrast {the consequences of different} learning {mechanisms}. In each of these categories, we have described in detail the prototypical models {and the} connections and departures between each. From these formulations, we have provided some of the important insights gained through studying these models. This includes some of the technical details concerning the development and analysis of these new models, as well as the current direction of study within particular classes of models. 

{Important too is the distinction and connection between mathematical and biological insights. Hence our effort to keep all technical mathematical questions raised connected to the ecological systems that {initiate our} motivation. This ranges from mathematical questions of existence and well-posedness (Sections \ref{sec:RoT} and \ref{sec:wellposedness}) to questions of more direct biological consequence (Sections \ref{sec:measuresofsuccess}-\ref{sec:critdom}).} To help provide direction for future study, we propose a wide variety of new models and related extensions to existing models throughout the manuscript. {To motivate researchers currently in the field or in adjacent fields, we propose many open problems related to all of the content explored. This includes specific, technical questions, but also includes a large number of open-ended questions that do not necessarily have a clear \textit{yes} or \textit{no} answer. These problems highlight how primitive some of the existing results are when compared to more mature areas of study, while emphasizing how much room there is for growth {and further study}.}

While we do not assert that we have provided a complete description of all existing cognitive mechanisms included within a diffusion-advection equation framework, we have made a substantial effort to include {a majority} of the common tools and techniques used. In cases where the big picture is perhaps treated as more important than the fine-grain detail, we have provided {many} relevant reference {materials} to accommodate further reading.

{From} the points raised above and the precise details found throughout this review, it is clear that these {knowledge-based movement} models and their extensions {will} have a broad impact for applied mathematicians and biologists alike. {M}any insights {are} provided here, but there are many more connections to be made. This includes a rich, diverse, and challenging branch of mathematical models which will require study from many different perspectives. New insights can be found through a more detailed exploration using existing mathematical techniques, while further insights will require {the development of} {more novel} tools and techniques, leaving much room for up-and-coming researchers to become pioneers in this growing field of study. {Complimentary to this}, mathematical explorations will be aided greatly by the contribution of knowledgeable biologists who can help make connections between analytical insights and biological ones, help make models biologically reasonable while favouring simplicity, and aid in the further development of new models and extensions beyond those discussed here. We hope this review will encourage new researchers to contribute to this exciting new {intersection} of mathematical biology and partial differential equations.

\section*{Acknowledgements}

We thank our reviewers for their time and effort in providing insightful suggestions {and critical feedback}, which greatly improved the quality of this manuscript. We also thank our colleagues Peter Thompson and Xiunan Wang for their constructive discussions in creation of Figure \ref{ConDia} and ultimately catalyzing our motivation to produce an independent document focusing on reaction-advection-diffusion equations. A special thank you to Sophia Salmaniw for painting the graphics used in Figure \ref{ConDia}.

Finally, we warmly thank Mark Lewis for his thoughtful feedback and for sharing his wealth of knowledge in the field of mathematical biology. We dedicate this work in honour of his 60th birthday.

Hao Wang's research was partially supported by NSERC Individual Discovery Grant RGPIN-2020-03911 and Discovery Accelerator Supplement Award RGPAS-2020-00090. Yurij Salmaniw was partially supported by NSERC Scholarship PGSD3-535063-2019, President's Doctoral Prize of Distinction, Josephine Mitchell Research Prize, and multiple Alberta Graduate Excellence Scholarships.

\begin{appendices}

\section{Equivalence of Models}\label{sec:equivalence}

In some cases, the models presented in this review can be reformulated into an equivalent model. We first write the full problem studied in \cite{Shi_JDDE}:
\eql{
\begin{cases}
\frac{\p u}{\p t} = d_1 \D u + d_2 \grad \cdot ( u \grad v) + f(u), \quad x \in \O, \ t>0, \cr
\frac{\p u}{\p \mathbf{n}} = 0 , \hspace{4.4cm} x \in \p \O, \ t>0, \cr
u(x,t) = \eta(x,t), \hspace{3.1cm} x \in \O, \ -\i < t \leq 0,
\end{cases}
}{Shi-fullmodel}
in a smooth, bounded domain $\O$, $\eta(x,t)$ is given initial data, and $v(x,t)$ is defined as
\eq{
v(x,t) = \mathcal{G} ** u = \int_{- \i} ^t \int_\O G(d_3, x,y,t-s) \mathcal{G}(t-s) u(y,s) dy ds ,
}
where $G$ is the Green's function for the heat equation in $\O$ subject to homogeneous Neumann boundary data, {$d_3$ is the diffusion rate for the Green's function,} and $\mathcal{G}$ is either the weak or strong kernel defined in \eqref{weakstrongkernel}. We first state Lemma 2.1 found in \cite{Shi_JMB}. The lemma is stated as follows.
\begin{lem}
Suppose that kernel $\mathcal{G}_w(t)$ is chosen to be the weak kernel defined in \eqref{weakstrongkernel} and define
\eq{
v(x,t) = \mathcal{G}_w ** u(x,t) = \int_{-\i}^t \int_\O G(d_3,x,y,t-s) \mathcal{G}_w (t-s) u(y,s) dy ds,
}
where $G$ is the Green's function for the heat equation subject to homogeneous Neumann boundary data. Then,
\begin{enumerate}
\item If $u(x,t)$ is the solution of \eqref{Shi-fullmodel}, then $(u(x,t), v(x,t))$ is the solution of
\eql{
\begin{cases}
\frac{\p u}{\p t} = d_1 \D u + d_2 \grad \cdot( u \grad v) + f(u), \hspace{3.3cm} x \in \O, \ t>0, \cr
\frac{\p v}{\p t} = d_3 \D v + \tau^{-1} ( u - v), \hspace{4.6cm} x \in \O,\ t>0, \cr
\frac{\p u}{\p \mathbf{n}} = \frac{\p v}{\p \mathbf{n}} = 0 , \hspace{6.4cm} x \in \p \O, \ t>0, \cr
u(x,t) = \eta (x,t), \hspace{6cm} x \in \O, \ t \leq 0, \cr
v(x,0) = \tau^{-1} \int_{-\i}^0 \int_\O G(x,y,-s) e^{s \tau^{-1}} \eta(y,s) dy ds,  \hspace{1.0cm} x \in \O, \ t \leq 0.
\end{cases}
}{Shi-fullmodel2}
\item If $(u(x,t), v(x,t))$ is a solution of
\eql{
\begin{cases}
\frac{\p u}{\p t} = d_1 \D u + d_2 \grad \cdot ( u \grad v) + f(u), \quad\quad x \in \O,\ t\in \mathbb{R}, \cr
\frac{\p v}{\p t} = d_3 \D v + \tau^{-1}(u-v), \hspace{2.1cm} x \in \O, \ t\in \mathbb{R}, \cr
\frac{\p u}{\p \mathbf{n}} = \frac{\p v}{\p \mathbf{n}} = 0 , \hspace{3.9cm} x \in \p \O, \ t\in \mathbb{R},
\end{cases}
}{Shi-fullmodel3}
then $u(x,t)$ satisfies equation \eqref{Shi-fullmodel} such that $\eta (x,s) = u(x,s)$, $-\i < s < 0$. In particular, if $(u(x), v(x))$ is a steady state of \eqref{Shi-fullmodel3}, then $u(x)$ is a steady state of \eqref{Shi-fullmodel}; if $(u(x,t), v(x,t))$ is a periodic solution of \eqref{Shi-fullmodel3}, then $u(x,t)$ is a periodic solution of \eqref{Shi-fullmodel}. 
\end{enumerate}
\end{lem}

This is an interesting result for two reasons. First, it is interesting to see that model \eqref{Shi-fullmodel} is actually equivalent to a Keller-Segel chemotaxis model when the weak kernel is chosen. Second, as a result of this first fact, one can then use the huge body of literature studying chemotaxis models in order to gain insights into this new delay partial differential equation. In the case where one chooses the strong kernel, there is another equivalent system consisting of $3$ equations and similar insights can be gathered. This is Lemma 2.2 in \cite{Shi_JMB}, which we do not provide here.

\end{appendices}

\bibliographystyle{acm}
\bibliography{references}
\end{document}